\documentclass[10pt]{amsart}
\usepackage{amssymb,latexsym,amsmath,psfig,amscd,epsfig,mathrsfs,yfonts}
\input xy
\xyoption{all}

\newcommand{\ZZ}{\mathbb Z}
\newcommand{\QQ}{\mathbb Q}
\newcommand{\CC}{\mathbb C}

\newcommand{\N}{\mathbb N}

\newcommand{\proj}{\mathbb P}
\newcommand{\grot}{\mathcal{A}}

\newcommand{\spc}{\vspace{3mm}}
\newcommand{\Sheaf}{\mathscr}
\newcommand{\beq}{\begin{eqnarray}}
\newcommand{\beqn}{\begin{eqnarray*}}
\newcommand{\eeq}{\end{eqnarray}}
\newcommand{\eeqn}{\end{eqnarray*}}
\newcommand{\nc}{non-commutative \hspace{0.5mm}}
\newcommand{\NC}{Non-commutative \hspace{0.5mm}}            
\newcommand{\QED}{\begin{flushright}$\blacksquare$\end{flushright}}
\newcommand{\ul}{\underline}

\newtheorem{thm}{Theorem}[section]

\newtheorem{lem}[thm]{Lemma}
\newtheorem{prop}[thm]{Proposition}
\newtheorem{cor}[thm]{Corollary}
\newtheorem{ex}{Example}
\newtheorem{defn}[thm]{Definition}
\theoremstyle{remark}
\newtheorem{rem}[thm]{Remark}

\newcommand{\sheaf}{\mathcal{O}}
\newcommand{\map}{\longrightarrow}

\newcommand{\lk}{{}_{A}k}

\begin{document}

\title[header at the top of the page]{On some approaches towards Non-commutative Algebraic Geometry}
\author{Snigdhayan Mahanta}
\address{Max-Planck-Institut f\"ur Mathematik, Bonn, Germany.}
\email{smahanta@mpim-bonn.mpg.de}

\begin{abstract}
The works of R. Descartes, I. M. Gelfand and A. Grothendieck have convinced us that commutative rings should be thought of as rings of functions on some appropriate (commutative) spaces. If we try to push this notion forward we reach the realm of {\it \NC\;Geometry}. The confluence of ideas comes here mainly from three seemingly disparate sources, namely, quantum physics, operator algebras (Connes-style) and algebraic geometry. Following the title of the article, an effort has been made to provide an overview of the third point of view. Since na\" ive efforts to generalize commutative algebraic geometry fail, one goes to the root of the problem and tries to work things out ``categorically''. This makes the approach a little bit abstract but not abstruse. However, an honest confession must be made at the outset - this write-up is very far from being definitive; hopefully it will provide glimpses of some interesting developments at least.
\end{abstract}

\thanks{This is a ``Diplom'' thesis carried out under the patronage of Max-Planck-Institut f\" ur Mathematik, Bonn. It is just an informal cruise through the subject and no originality, in any form, is claimed.}%

\maketitle

\begin{center}
\today
\end{center}



\section{Introduction}

\vspace{5mm}

It behooves me to explore a little bit of history behind the
  subject . However, it will be short as this is not meant to be a discourse on
  mathematical history and, in the process, significant contributions of several
  mathematicians spanning more than two decades of mathematical research will
  have to be overlooked. Sincere apologies are offered to them. 

 In a seminal paper entitled {\it
  Faisceaux alg\'ebriques coh\'erents} \cite{Ser} in 1955 {\it Serre} had
introduced the notion of {\it coherent} sheaves. Hindsight tells us that the
seeds of \nc algebraic geometry were sown in this work, though it is not clear
if Serre actually had such an application in his mind. In fact, it was {\it
  Manin} around 1988 \cite{Man2}, who first propounded the idea of dropping the
commutativity hypothesis. {\it Pierre Gabriel} in his celebrated work on
abelian categories \cite{Gab} in 1962 had also proved some reconstruction theorems
for noetherian schemes, which were general enough to have an appeal to the
mathematicians looking to go beyond the commutative framework. Apparently, the
  real upsurge in such activities started taking place roughly 20 years back,
  with the works of {\it Artin} and {\it Schelter} on the classification of
  \nc rings into some non-trivial classes. Of course, it was not entirely an
  aimless pursuit as quantum physics had been a steady source of problems of
  \nc nature. There seems to be a plethora of proposals on the way to conduct
  research in this area. One of the most successful approaches was provided by
  {\it Connes} (see \cite{Con} and \cite{Man1} for a broader perspective) but that will be outside the
  purview of our discussion as we have an algebraic epithet in our title. Some
  more historical anecdotes and trends in current research can be found in the
  last chapter. Now we will take a brief look at the contents of each chapter. 

\vspace{1mm}
In chapter 1 we take a look at the description of {\it Artin-Schelter regular
  algebras} of dimension $3$ after {\it Artin, Tate} and {\it Van den
  Bergh}. We see how they arise from triples $(E,\sigma,\Sheaf{L})$, which are
  well studied objects in commutative algebraic geometry.

In chapter 2 we almost entirely dedicate ourselves to {\it Rosenberg's}
reconstruction of varieties from their categories of quasi-coherent
sheaves. In this, rather arid, excerpt we familiarize ourselves with {\it Grothendieck Categories} and find out how one can associate a ringed space to such a
category. 

In chapter 3 an effort has been made to understand the preliminaries of the
model of {\it non-commutative projective geometry} spearheaded by {\it Artin}
and {\it Zhang}. Towards the end we also include a discussion on \nc proper
varieties which is a cumulative effort of works done by many different
mathematicians. 

In chapter 4 we provide, a rather sketchy, outline of a model of \nc algebraic
geometry initiated by {\it Laudal}. It has the advantage of having a good
theoretical basis for working infinitesimally but it does not seem as natural as
the model described in chapter 3. Unfortunately we cannot even describe projective geometric objects in this approach as getting to the
affine ones is quite an uphill task in itself. 

In chapter 5 we shall catch a glimpse of some other approaches in this
direction which have had varied levels of success.

\vspace{1mm}
Finally, it must be bourne in mind that it is a brief survey of \nc
algebraic geometry and hence, as asserted earlier, very far from being
exhaustive. References to the suitable articles have been provided amply and
their ubiquity, coupled with the general expository nature of this write-up, account for the lack of proofs. One final remark: this
tome (except for chapter 1) is redolent of or stinks of categories as one's tastes might be.  


\newpage

In a paper entitled {\it Some Algebras Associated to Automorphisms of Elliptic Curves} Artin, Tate and Van den Bergh \cite{ATV1} had given a nice description of \nc algebras which should, in principle, be algebras of functions of some nonsingular ``\nc schemes''. In the commutative case, nonsingularity is reflected in the regularity of the ring. However, this notion is insufficient for \nc purposes. So Artin and Schelter had given a stronger regularity condition which we call {\it Artin-Schelter}(AS)-{\it regularity} condition. The main result of the above mentioned paper says that AS-regular algebras of dimension 3 (global dimension) can be described neatly as some algebras associated to automorphisms of projective schemes, mostly elliptic curves. And also such algebras are both left and right Noetherian. This sub-section is entirely based on the contents of \cite{ATV1}.
\vspace{2mm}

To begin with, we fix an algebraically closed field $k$ of
characteristic $0$. We shall mostly be concerned with $\N$-graded
$k$-algebras $A\;=\;\underset{i\geqslant 0}{\oplus} A_i$, that are
finitely generated in degree $1$, with $A_0$ finite dimensional as a
$k$-vector space. Such algebras are called {\it finitely graded} for
short, though the term could be a bit misleading at the first sight.
A finitely graded ring is called {\it connected graded} if
$A_0\;=\;k$. $A_{+}$ stands for the two-sided {\it augmentation
ideal} $\underset{i>0}{\oplus} A_i$.

\spc

\begin{defn} (AS-regular algebra)

A connected graded ring $A$ is called Artin-Schelter (AS) regular of dimension $d$ if it satisfies the following conditions:

(1) $A$ has global dimension $d$.

(2) GKdim($A$) $<$ $\infty$.

(3) $A$ is AS-Gorenstein.

\end{defn}

It is worthwhile to say a few words about {\it Gelfand-Kirillov\;dimension} (GKdim) and the {\it AS-Gorenstein} condition of algebras.

\spc

Take any connected graded $k$-algebra $A$ and choose a finite dimensional $k$-vector space $V$ such that $k[V]\;=\;A$. Now set $F^{n}A\;=\;k\;+\;\sum_{i=1}^{n} V^i$ for $n \geqslant 1$. This defines a filtration of $A$. Then the $GKdim(A)$ is defined to be

\begin{eqnarray*}
GKdim(A)\;=\;\underset{n}{lim sup}\frac{ln(dim_k F^{n}A)}{ln(n)}.
\end{eqnarray*}

Of course, one has to check that the definition does not depend on the choice of $V$.

\begin{rem}
{\it Bergman} \cite{KL} has shown that $GKdim$ can take any real number $\alpha \geqslant 2$. However, if $GKdim\;\leqslant\;2$, then it is either $0$ or $1$.
\end{rem}

There are some equivalent formulations of the {\it AS-Gorenstein} condition available in literature. We would just be content by saying the following:

\newpage

\begin{defn} (AS-Gorenstein condition) \label{Gorenstein}

A connected graded $k$-algebra $A$ of global dimension $d<\infty$ is AS-Gorenstein if 
\begin{eqnarray*}
Ext^i_A(k,A)\;=\;0\;for\;i\neq d\;and\;Ext^d_A(k,A)\simeq k
\end{eqnarray*}
\end{defn}

All regular commutative rings are AS-Gorenstein, which supports our conviction that the AS-Gorenstein hypothesis is desirable for \nc analogues of regular commutative rings. Further, note that the usual {\it Gorenstein} condition (for commutative rings) requires that they be Noetherian of finite injective dimension as modules of themselves.

Now we take up the task of describing the minimal projective resolution $(\star)$ of an AS-regular algebra of dimension $d\;=\;3$. As a fact, let us also mention that the global dimension of a graded algebra is equal to the projective dimension of the left module $\lk$. Let

\beq
0 \map P^d\map\dots\overset{f_2}{\map}P^1\overset{f_1}{\map}P^0\map \lk\map 0
\eeq

be a minimal projective resolution of the left module $\lk$. $P^0$ turns out to be $A$; $P^1$ and $P^2$ need a look into the structure of $A$ for their descriptions. Suppose $A\;=\;T/I$, where $T\;=\;k\{x_1,\dots,x_n\}$ is a free associative algebra generated by homogeneous elements $x_i$ with degrees $l_{1j}$ (also assume that $\{x_1,\dots,x_n\}$ is a minimal set of generators). Then

\beq
P^1\;\approx\;\underset{j=1}{\overset{n}{\oplus}}\;A(-l_{1j})
\eeq

The map $P^1\map P^0$, denoted $x$, is given by right multiplication with the column vector $(x_1,\dots,x_n)^t$.

Coming to $P^2$, let $\{f_j\}$ be a minimal set of homogeneous generators for the ideal $I$ such that $deg\;f_j\;=\;l_{2j}$. In $T$, write each $f_j$ as

\beq
f_j\;=\;\underset{j}{\sum}m_{ij}x_j
\eeq

where $m_{ij}\;\in\;T_{l_{2i}-l_{1j}}$. Let $M$ be the image in $A$ of the matrix $(m_{ij})$. Then

\beq
P^2\;\approx\;\underset{j}{\oplus} A(-l_{2j})
\eeq

and the map $P^2\map P^1$, denoted $M$, is just right multiplication by the matrix $M$.

\spc

In general, it is not so easy to interpret all the terms of the resolution $(1)$. However, for a regular algebra of dimension $3$, the resolution looks like

\beq
0 \map A(-s-1)\stackrel{x^t}{\map} A(-s)^r\stackrel{M}{\map} A(-1)^r\stackrel{x}{\map} A\map \lk \map 0 \hspace{10mm}
\eeq

where $(r,s)\;=\;(3,2)$ or $(2,3)$.  {\it Thus such an algebra has $r$ generators and $r$ relations each of degree $s$, $r+s\;=\;5$}. Set $g\;=\;(x^t)M$; then

\beq
g^t\;=\;((x^t)M)^t\;=\;QMx\;=\;Qf \hspace{10mm}
\eeq

for some $Q \in GL_r(k)$.

Now, with some foresight, we introduce a new definition, that of a {\it standard algebra}, in which we extract all the essential properties of AS-regular algebras of dimension $3$.

\begin{defn} An algebra $A$ is called \underline{standard} if it can be presented by $r$ generators $x_j$ of degree 1 and $r$ relations $f_i$ of degree s, such that, with M defined by $(3)$, $(r,s)\;=\;(2,3)$ or $(3,2)$ as above, and there is an element $Q \in GL_r(k)$ such that $(6)$ holds.
\end{defn}

\begin{rem}
For a standard algebra $A$, $(5)$ is just a complex and if it is a resolution, then $A$ is a regular algebra of dimension 3.
\end{rem}

\begin{center}
{\bf Multilinearization and a Moduli Problem}
\end{center}

\spc

We have some algebraic objects, namely \nc $k$-algebras, and we want to associate some geometric data to them. Ideally we would like to think of the algebra as the ``coordinate ring'' of some ``space'', but this is a utopia at the moment. Nevertheless, we can associate some other geometric objects to them, namely a family of schemes $\{\Gamma_d\}$, which is called the {\it multilinearization} of the algebra. Now we shall try to understand the interplay between the two. 

\begin{eqnarray*}
\N-graded\;\nc\;k-algebras \longleftrightarrow families\;of\;schemes\;\{\Gamma_d\}
\end{eqnarray*}

\spc

However, the figure above is not as nice as it looks; the arrow is not really reversible. If one starts with an algebra $A$ and passes on to its multilinearization and then tries to recover the algebra from it, one gets another algebra $B$ together with an algebra map $A\;\map\;B$ which is bijective in degree 1. 

\spc
\noindent
\underline{\bf Algebra to Family of Schemes}
\vspace{1mm}

Let $A=T/I$ be a graded k-algebra, where $T=k(x_0,\dots,x_n)$ is the free associative algebra generated by $x_0,\dots,x_n$ of degree 1 and $I$ is a homogeneous ideal of $T$.

Let $V=T_1$ and set $P=\proj(V)=\proj^n_k$. Let $\sheaf(1,\dots,1)=pr_1^\ast\sheaf_P(1)\otimes\dots\otimes pr_d^\ast\sheaf_P(1)$ denote the twisting sheaf on the product (fibre product over $Spec\, k$) of $d$ copies of $P$. Every homogeneous element $f\in T_d = T_1^{\otimes d}$ defines a global section $\tilde{f}$ of this sheaf. Note that $A$ was defined as a quotient of $T$ by a homogeneous ideal $I$. Let $\tilde{I}_d = \{\tilde{f}:f\in I_d\}$. Then we define 

\beq
\Gamma_d = \mathcal{Z}(\tilde{I}_d) \subset P^d
\eeq

\noindent
where $\mathcal{Z}$ denotes the scheme of zeros. Of course, with natural conventions we have,

\beq
\tilde{(fg)}(v_1,\dots,v_p,w_1,\dots,w_q)=\tilde{f}(v_1,\dots,v_p)\tilde{g}(w_1,\dots,w_q)
\eeq

\noindent
for $f\in T_p,\;g\in T_q$  

\spc

It must be mentioned here that actually the family $\tilde{I}_d$ is called the {\it multilinearization} of $A$ as it consists of multilinear forms on $V$ but by abuse of notation, which is very much in vogue, we call their scheme of zeros multilinearization as well. Let us rush through some of the properties of this family of schemes so defined. Rather formal arguments will lead to a proof of most of them \cite{ATV1} and hence will not be provided here. 

\begin{prop}

1. For any $d$, $\Gamma_{d+1} \subset (P\times\Gamma_d) \cap (\Gamma_d\times P)$\hspace{1mm} and equality holds if $I_{d+1}=T_1I_d+I_dT_1$ (for instance if $I$ is generated in degrees $\leqslant d$).

\vspace{1mm}
2. Let $P_i$ denote the $i-th$ factor of the product $P^d$ and for $1\leqslant i < j\leqslant d$, let $pr_{ij} = pr^{(d)}_{ij}$ denote its projection to the product $\underset{\nu=i}{\overset{j}{\prod}}P_\nu$. Then $pr_{ij}(\Gamma_d)$ is a closed subset of $\Gamma_{j-i+1}$.

\end{prop}

\begin{prop}  \label{prop_def_pi}

Let $0\leqslant i \leqslant d$ and let $\pi_i : \Gamma_{d+1}\map\Gamma_d$ denote the projection (dropping the $(i+1)-th$ component)

\beq
(p_1,\dots,p_i,p_{i+1},p_{i+2},\dots,p_{d+1}) \longmapsto (p_1,\dots,p_i,p_{i+2},\dots,p_{d+1})
\eeq

\spc

1. The fibres of $\pi_i$ are linear subspaces of $P$.

2. Let $p\in\Gamma_d$, and let $L$ be the fibre of $\pi_i$ at $p$. If $dim L \leqslant 0$ (this encompasses the fact that the fibre at a point could be empty), then $\pi_i$ is a closed immersion locally in a neighbourhoud of $p\in\Gamma_d$.
\end{prop}

\begin{prop} \label{semi_standard}

1. Assume that for some $d$, $pr^{(d+1)}_{1d}$ defines a closed immersion from $\Gamma_{d+1}$ to $\Gamma_d$, thus identifying $\Gamma_{d+1}$ with a closed subscheme $E\subset\Gamma_d$. Then $\Gamma_{d+1}$ defines a map $\sigma:E \map \Gamma_d$, by the rule

\beq
\sigma(p_1,\dots,p_d) = (p_2,\dots,p_{d+1})
\eeq

\spc
where $(p_1,\dots,p_{d+1})$ is the unique point of $\Gamma_{d+1}$ lying above $(p_1,\dots,p_d)\in\Gamma_d$.

2. If, in addition to this, $\sigma(E)\subset E$ and if $I$ is generated in degree $\leqslant d$, then $pr^{(n)}_{ij}:\Gamma_n \map E$ is an isomorphism for every $n \geqslant d+1$.

\end{prop}

\spc

Now we introduce the definition of {\it point modules}.

\begin{defn}

A graded right $A$-module $M$ will be called a {\it point module} if it satisfies the following conditions:

1. $M$ is generated in degree zero.

2. $M_0 = k$, and

3. $dim_k M_i = 1$ for all $i \geqslant 0$. 

By adding one more adjective, namely truncated, we arrive at the notion of a \underline{truncated point module}. A truncated point module of length $d+1$ is a point module whose Hilbert function is truncated, {\it i.e.,}

\begin{eqnarray*}
dim M_i = 
    \begin{cases}
      1&   \text{if $0\leqslant i \leqslant d$,} \\
      0&   \text{otherwise.}
    \end{cases}
\end{eqnarray*}

 Put differently, a {\it point module} can be thought of as a module, generated in degree zero, whose Hilbert series is $(1-t)^{-1}$ (and for a truncated point module it is just the polynomial $1+t+\dots +t^d$). 

\end{defn}

Our next aim is to unravel the following statement: 

\vspace{2mm}
\noindent
{\bf $\Gamma_d$ represents the functor of flat families of truncated point modules of length $d+1$ parametrized by $Spec\,R$.}

\spc

As a digression let me say a few words about {\it moduli problems}. For reasons not very hard to see, it is desirable to come up with a ``space'' which parametrizes certain families of objects defined over schemes modulo some equivalence relation. Let $k$ be a field and let $Sch_{k}$ denote the category of schemes over $k$ (or $Spec\;k$ if you like). Let $F$ be a functor from $Sch_{k}^o$ to $Sets$ (category of sets), 

\beq 
\begin{split}
F:Sch_{k}^o &\map Sets \\
      B     &\longmapsto  F(B)
\end{split}
\eeq

where $F(B)$ is the set of equivalence classes of families of objects parametrized by the scheme $B$. Such an $F$ is called a moduli problem. A typical example will be the problem where we associate to a scheme the family of elliptic curves with an {\it n-level structure} modulo isomorphism or the {\it Picard functor} for a curve. 

\spc

An ideal solution to a moduli problem would be a natural isomorphism $\Phi$ to the functor of points of a scheme, say $M$. If such an $M$ exists, then it is called the {\it fine moduli space} of the corresponding moduli problem. Let me remind you that a scheme $X$ is no more than its functor of points $h_X\;=\;Hom({}_{-},X)$ by Yoneda's Lemma. But we are not always fortunate enough to be able to construct a fine moduli space for a moduli problem and so we define what is called a {\it coarse moduli space}. However, I am not going to go into these discussions over here. 

\spc 

It should also be bourne in mind that if $X$ is a $k$-scheme, then $Hom(Spec\,R,X)$, denoted $X(R)$, is called the set of $R$-valued points of the scheme $X$, where $R$ is a commutative $k$-algebra.

\spc

Having gone through the general set-up the first thing that we need to do is to formulate a notion of a family of point modules. We shall work over the category of affine $k$-schemes, which is equivalent to the category of commutative $k$-algebras, denoted $Alg_k$. 

We once again turn back to our graded $k$-algebra $A$ (not necessarily commutative). Let $R\in Alg_k$ and consider a graded $R\otimes_k A$-module $M$ generated in degree 0 and satisfying $M_0 = R$. By restriction of scalars via the canonical map $R \longrightarrow R\otimes_k A$, we consider $M$ as an $R$-module. Since we are trying to define a family of point modules, it is expected that locally {\it i.e.,} for all $p \in Spec\, R$, $M_p$ behave like a point module. Hence, we also ask that for every degree $d$, $\tilde{M}_d$ (the sheaf associated to the graded piece $M_d$ of degree $d$) be a locally free sheaf of rank 1 over $Spec\;R$. 

To define a family of truncated point modules of length $d+1$ we just ask that $M_j$ be locally free of rank 1 for all $0\leqslant j\leqslant d$ and otherwise the stalk be 0 in some sense.  

To phrase it in the moduli problem terms we say that our functor $F$ associates to any affine scheme over $k$, $Spec\;R$, the isomorphism classes of all such $R\otimes_k A$-modules $M$ (the truncated version). One has to take the obvious notion of isomorphism amongst such modules. 

At the end of all these discussions we may conclude that the statement in bold letters is equivalent to saying, 

\begin{prop}

$\Gamma_d$ is the fine moduli space for the moduli problem $F$.

\end{prop}

\noindent
{\bf Sketch of the Proof:}
Our aim is to show that the $R$-valued points of $\Gamma_d$ correspond to the isomorphism classes of flat families of truncated point modules of length $d+1$, parametrised by $Spec\,R$.

We start with a family of truncated point modules of length $d+1$ parametrised by $Spec\,R$. So, by definition this is a graded $R\otimes A$-module $M$, generated in degree zero, such that $M_0=R$ and that $M_i$ is locally free of rank 1 for each $i,1\leqslant i \leqslant d$. Now choose a covering $Spec\, R = \underset{j}{\cup} Spec\, R_{f_j}$, such that each $(M_{f_j})_i$ is a free $R_{f_j}$-module of rank $1$. To simplify notations, let me assert right here that if we can find $Spec\, R_{f_j} \map \Gamma_d$ for each $j$, then they will glue properly to give rise to a map $Spec\, R \map \Gamma_d$ and hence, an $R$-valued point of $\Gamma_d$. So we might as well pretend that $R$ were equal to one of its localizations and hence, ``drop the subscripts $f_j$''. 

We have a collection of free $R$-modules of rank 1; we promptly choose a basis for each one of them. Let $m_i$ be the chosen basis element for $M_i$. If you have all our hypotheses at your fingertips, you can easily recall that $A$ was generated by $x_j,\;0\leqslant j \leqslant n$, which were elements of degree 1. Now we write out the products of the bases $m_i$ by $x_j$.

\beq
m_{i-1}x_j\,=\,m_i a_{ij}
\eeq

\noindent
for some $a_{ij}\in R$. In this way we obtain a set of $d$ points 

\beq
a_i\,=\,(a_{i0},\dots,a_{in})\in P \;\forall\;0\leqslant i \leqslant d
\eeq

\noindent
Evidently $a=(a_1,\dots,a_d)$ gives an $R$-valued point of $P^d$. But I want to claim that actually $a\in\Gamma_d(R)$. Take any $a\in P^d(R)=Hom(Spec\,R,P^d)\simeq Hom(T_d,R)$ obtained in the manner described above. Let me remind you that $T_d=T_1^{\otimes d}$ corresponds to the scheme $P^d$. Now for any element $f\in R\otimes T$ of degree $d$, simply by looking at the definition of $a_i$ we can convince ourselves that

\beq
m_0 f\,=\,m_d\tilde{f}(a) \;\;\; \label{module_defn}
\eeq

But if $f\in R\otimes I_d$ then $f=0$ in $R\otimes A$ and hence $\tilde{f}(a)=0$. And $\Gamma_d$ was defined precisely as the scheme of zeros of $\tilde{I_d}$. So $a$ is indeed an $R$-valued point of $\Gamma_d$. 

\vspace{2mm}
\noindent
[If the symbols confuse you, let me put it for you differently. The element $a$ should be a map from $T\map R$. Take any $f\in I_d\subset T_d$. It is nothing but a polynomial in $x_0,\dots,x_n$ with coefficients in $k$, which we denote $\tilde{f}$. Tensor it with $R$ to obtain a polynomial with coefficients in $R$. Now ``evaluate'' it at $a$. This evaluation map has a kernel and that is precisely $I_d$ (or $\tilde{I_d}$ if you like).] 

\spc

Now the passage from $\Gamma_d$ to $M$ or schemes to modules is just a matter of walking back. We start with a point of $\Gamma_d(R)$ {\it i.e.,} $\rho\, :\, Spec\, R \map \Gamma_d$. Then consider the following sheaf.

\begin{eqnarray*}
\Sheaf{L}_j = i^\ast\sheaf (\underset{\text{$j$ times}}{\underbrace{1,\dots,1}},0,\dots,0)
\end{eqnarray*}

where $i:\Gamma_d\hookrightarrow P^d$, and put $\Sheaf{L}=\underset{j=0}{\overset{d}{\oplus}}\Sheaf{L}_j$. 

We just pull back this sheaf $\Sheaf{L}$ on $\Gamma_d$ via $\rho$ to obtain the desired truncated point module. Imitating the authors, I allow myself the luxury of cutting some corners and summing it up by asserting that the two ``functors'' defined above are ``quasi-inverses'' of each other.

\QED

\spc

\begin{rem}

This functor of points description of the family of schemes $\{\Gamma_d\}$ says that it is intrinsic to the algebra $A$, and does not depend on the presentation of $A$ as a quotient of $T$.

The family $\{\Gamma_d\}$ along with the maps $\pi$ defined in [\ref{prop_def_pi}] forms an inverse system and the inverse limit gives us an object whose ``points'' correspond to ``point modules''. 

\end{rem}

\noindent
\underline{\bf From Family of Schemes back to Algebra}

\spc

The most na\" ive attempt to recover an algebra $A$ from its multilinearisation $\{\Gamma_d\}$ yields an algebra $B$ and a canonical morphism $A\map B$. In fact, we can associate canonically a graded algebra $B$ to any sequence of subschemes $\{Z_d\subset P^d\}$ having the property

\beq
pr_{1,d-1}(Z_d)\subset Z_{d-1}\;\;and\;\;pr_{2,d}(Z_d)\subset Z_{d-1}\;\forall \,d
\eeq

\spc

Let $L_d = \sheaf(1,\dots,1)\otimes_{\sheaf_{P^d}}\sheaf_{Z_d}$. Then set $B_d=H^0(Z_d,L_d)$. Since $pr_{1,i}(Z_{i+j})\subset Z_i$ and $pr_{i+1,i+j}(Z_{i+j})\subset Z_j$, we have $Z_i \times Z_j \hookleftarrow Z_{i+j}$. Applying the $H^0$ functor will just reverse the arrows and define for us a multiplication map $B_i\times B_j \map B_{i+j}$ using the isomorphism $pr_{1,i}^\ast(L_i)\otimes_{\sheaf_{P^d}}pr_{i+1,i+j}^\ast(L_j)\map L_{i+j}$. This makes $B=\underset{d}{\oplus}B_d$ into a graded associative algebra. It must be mentioned that $B$ need not be generated in degree 1. If it is so, then the sequence $\{Z_d\}$ can be ``properly contained'' in the multilinearization of $B$, $\Gamma_d(B)$.

\begin{prop}

Let $A=T/I$ be a quotient of a free associative ring and let $\{\Gamma_d\}$ be the multilinearisation of $A$. Let $B$ be the algebra associated to it. Then, there is a canonical morphism, $\phi:A\map B$, which is bijective in degree 1.

\end{prop}

\noindent
{\bf Proof:}

The functorial maps

\beq
H^0(P^d,\sheaf(1,\dots,1)) \map H^0(\Gamma_d,L_d)
\eeq

define a morphism at the level of the graded pieces and hence a morphism from the free algebra $T$ to $B$. This is just the restriction of sections. Once again the definition of $\Gamma_d$ as the scheme of zeros of sections $\tilde{I_d}$ forces $I$ to be in the kernel of this morphism and hence this map factors through $A=T/I$. 

In degree 1, $\Gamma_1$ is the dual projective space of $A_1=T_1/I_1$ and hence its sheaf of global sections $H^0(\Gamma_1,L_1)=B_1 \simeq A_1$ [Lefschetz Hyperplane Theorem]. (It is dual because $V=T_1^\ast$ and $\Gamma_1=\proj((T_1/I_1)^\ast)\subset \proj(V)=P$). 
\QED

This process of passing onto the multilinearization is rather unwieldy. Fortunately, we don't need to consider the entire family of schemes. We can concentrate on a smaller piece of datum in the form of a triple $(E,\sigma,\Sheaf{L})$, which will be taken up now. The recipe to obtain the triple from an algebra (if at all possible) is explained in the ``main result'' subsection. 

\spc

\begin{center}
{\bf Twisted Homogeneous Coordinate Rings}
\end{center}

\spc

Here we provide a very general recipe of manufacturing interesting \nc rings out of a completely ``commutative geometric'' piece of datum, called an {\it abstract triple}, which turns out to be an isomorphism invariant for AS-regular algebras.  

\begin{defn}
An abstract triple $\mathcal{T}=(X,\sigma,\Sheaf{L})$ is an assortment of a projective scheme $X$, an automorphism $\sigma$ of $X$ and an invertible sheaf $\Sheaf{L}$ on $X$.
\end{defn}

It is time to construct the {\it Twisted Homogeneous Coordinate Ring} $B(\mathcal{T})$ out of an abstract triple. For each integer $n\geqslant 1$ set 

\beq
\Sheaf{L}_n\;=\;\Sheaf{L}\otimes\Sheaf{L}^\sigma\otimes\dots\otimes\Sheaf{L}^{\sigma^{n-1}}
\eeq

where $\Sheaf{L}^\sigma\; :=\;\sigma^{\ast}\Sheaf{L}$. The tensor products are taken over ${\sheaf}_X$ and we set $\Sheaf{L}_0\;=\;{\sheaf}_X$. As a graded vector space, $B(\mathcal{T})$ is defined as

\begin{eqnarray*}
B(\mathcal{T})\;=\;\underset{n\geqslant 0}{\oplus} H^0(X,\Sheaf{L}_n)
\end{eqnarray*}

For every pair of integers $m,n\geqslant 0$, there is a canonical isomorphism

\beqn
\Sheaf{L}_m\otimes_{{}_k} \Sheaf{L}_n^{\sigma^m} \map \Sheaf{L}_{m+n}
\eeqn

and hence defines a multiplication on $B(\mathcal{T})$

\beqn
H^0(X,\Sheaf{L}_m)\otimes_{{}_k} H^0(X,\Sheaf{L}_n)\map H^0(X,\Sheaf{L}_{m+n}).
\eeqn

\begin{ex}

Let us compute (more precisely, allude to the computation of) the twisted homogeneous coordinate ring in a very simple case. Let $\mathcal{T}=(\proj^1,\sheaf(1),\sigma)$, where $\sigma(a_0,a_1)=(q a_0,a_1)$ for some $q\in k^\ast$. 

\end{ex}

We choose a parameter $u$ for $\proj^1$, so that the standard affine open cover $U=\proj^1\setminus \infty$, $V=\proj^1\setminus 0$ has rings of regular functions $\sheaf(U)=k[u]$ and $\sheaf(V)=k[u^{-1}]$. Now we can identify $\sheaf(1)$ with the sheaf of functions on $\proj^1$ which have at most a simple pole at infinity; in other words, it is a sub-sheaf of $k(u)=k(\proj^1)$ generated by $\{1,u\}$. It is now easy to see that $H^0(X,\sheaf(n))$ is spanned by $\{1,u,\dots,u^n\}$ and that, as a graded vector space $B(\proj^1,id,\sheaf(1))=k\{x,y\}$ (the free algebra over $k$ generated by $x$ and $y$ and not the usual polynomial ring), where $x=1$ and $y=u$, thought of as elements of $B_1=H^0(X,\sheaf(1))$. It should be mentioned that $\sigma$ acts on the rational functions on the right as $f^\sigma(p)=f(\sigma(p))$ for any $f\in k(\proj^1)$ and $p\in \proj^1$. It is clear from our concrete presentation that $\sheaf(1)^\sigma\cong\sheaf(1)$. So as a graded vector space $B(\mathcal{T})\cong B(\proj^1,id,\sheaf(1)$. However, the multiplication is twisted.

\beqn
y.x\; =\; y\otimes x^\sigma\; =\; u\otimes 1^\sigma\; =\; u\otimes 1\; =\; u\in H^0(\proj^1,\sheaf(2)).
\eeqn

On the other hand,  

\beqn
x.y\; =\; x\otimes y^\sigma\; =\; 1\otimes u^\sigma\; =\; 1\otimes qu\; =\; qu\in H^0(\proj^1,\sheaf(2)).
\eeqn

So we find a relation between $x$ and $y$, namely, $x.y \,-\,qy.x\,=\,0$ and a
little bit of more work shows that this is the only relation. So the twisted
homogeneous coordinate ring associated to $\mathcal{T}$ is

\beqn
B(\mathcal{T})\;=\; k\{x,y\}/(x.y\,-\,qy.x).
\eeqn

\noindent
A good reference for a better understanding of these rings is \cite{AV}.

\spc

\begin{center}
{\bf The Main Result}
\end{center}

 The main result of the paper \cite{ATV1} asserts that the AS-regular algebras of dimension $3$ are precisely the {\it non-degenerate standard algebras}. Let me now devote some time to the italicized words, which we have not seen yet and hopefully the reader will be convinced that they are more tractable objects and give us a better understanding of the AS-regular algebras of dimension $3$. 

\spc

Let $A$ denote a $k$-algebra which is presented in the form $A=T/I$, where $T$ is a \nc polynomial ring with $r$ generators of degree $1$ and $I$ is an ideal generated by $r$ linearly independent relations of degree $s$. We assume, moreover, that $(r,s)\,=\,(2,3)$ or $(3,2)$ keeping the definition of a standard algebra in mind, which in turn is invoked to concentrated on AS-regular algebras of dimension $3$.

Let $f=(f_1,\dots,f_r)^t$ be the column vector of defining relations {\it i.e.,} of generators of $I$. There are uniquely defined matrices $M,N$ with entries in the tensor algebra $T$ such that 

\beq
\text{$f = Mx$ and $f^t = x^t N.$}
\eeq

Let \underline{$P=\proj^{r-1}$}. By $P^d$ we mean $\underset{d\;times}{\underbrace{P\times_k\dots\times_k P}}$. We consider the multilinearization $\{\Gamma_d \subset P^d\}$. It follows from [\ref{semi_standard}] that the whole sequence of schemes is determined by $\Gamma = \Gamma_s$. Set $\pi_1 = pr_{1,s-1}\mid_{{}_{{}_{\Gamma}}}$ and $\pi_2 = pr_{2,s}\mid_{{}_{{}_{\Gamma}}}$ (refer to [\ref{prop_def_pi}]). Now let $E = \pi_1(\Gamma)$ and $E' = \pi_2(\Gamma)$. By definition, $E$ is the locus of points $(p_1,\dots,p_{s-1})\in P^{s-1}$ for which there exists $p_s \in P$ such that $(p_1,\dots,p_{s-1},p_s) \in \Gamma$. So $E$ is the locus of zeros in $P^{s-1}$ of the multihomogeneous polynomial $det \tilde{M} = 0$. Similarly, $E'$ is the locus of $det \tilde{N} = 0$. Thus $E$ is either all of $P^{s-1}$ if $det \tilde{M} \equiv 0$, or is a Cartier divisor in $P^{s-1}$, and similarly for $E'$.

\begin{defn}

We call an algebra $A$ semi-standard if the schemes $E$ and $E'$ are equal in $P^{s-1}$ or equivalently if 

\beq
det \tilde{M}\, =\, c.det \tilde{N}
\eeq

for some $c \in k^{\ast}$.

\end{defn}

The constant $c$ is independent of the change of basis in $x$ and $f$, provided that $det \tilde{M}$ is not identically $0$. As one can see, this definition does not really need the assumption that $(r,s) = (2,3)$ or $(3,2)$. 

Now, let $A$ be a semi-standard algebra, so that $E = E'$. We may view $\Gamma$ as the graph of a correspondence $\sigma$ on $E$ via the closed immersion $(\pi_1,\pi_2): \Gamma \map E\times E$. It should be noted that $\pi_2 = \pi_1\circ \sigma$. 
 
\begin{defn}
We call a semi-standard algebra non-degenerate if the $sigma$ mentioned above is actually an automorphism of $E$, and degenerate otherwise.
\end{defn}

If there is any sanity in the nomenclature, then one should expect standard algebras to be semi-standard. This is indeed the case and its verification is just a matter of playing with the definitions. 

The good news is that we are already in sight of an abstract triple. If $A$ is a non-degenerate semi-standard algebra then we can associate the triple $\mathcal{T}(A) = (E,\sigma,\Sheaf{L})$, where $E=\pi_1(\Gamma)\stackrel{i}{\hookrightarrow} P^{s-1}$, $\sigma$ is the automorphism given by the non-degeneracy of the algebra and $\Sheaf{L}$ is contrived from the twisting sheaf as $\psi^{\ast}\sheaf_{P}(1)$, where $\psi:E \map P\, =\,\mathbb{P}^{r-1}$ is the morphism defined by 

\beq
\text{}
\psi =
    \begin{cases}
     i:E \map \mathbb{P}^2 &  \text{if $(r,s) = (3,2),$} \\
     pr_1\circ i: E \map \mathbb{P}^1 & \text{if $(r,s) = (2,3).$}
    \end{cases}
\text{}
\eeq

It must be bourne in mind that $(r,s)\, =\, (2,3)$ or $(3,2)$ and hence, $P^{s-1}\, = \mathbb{P}^1\times_k \mathbb{P}^1$ or $\mathbb{P}^2$ respectively. Hence, $\psi$ is determined by the sections of $\Sheaf{L}$. This procedure of obtaining abstract triples can also be applied to any family of schemes $\{\Gamma_d \subset P^d\}$, satisfying $\Gamma_d \cong \Gamma_s$ $\forall\; d\geqslant s$ and $\pi_1(\Gamma_s) = \pi_2(\Gamma_s) \subset P^{s-1}$.

It follows from the definition that AS-regular algebras of dimension $3$ are standard. And it has also been shown in \cite{ATV1} that AS-regular algebras are non-degenerate. So we can associate a triple to such an algebra. The regularity of the algebra is reflected in the following property [\ref{regularity}] of the triple, which is also called the {\it regularity} property of the triple for obvious reasons. 

\beq \label{regularity}
\text{}
  \begin{cases}
  (\sigma - 1)^2 [\Sheaf{L}]\, =\, 0&  \text{if $r = 3$} \\
  (\sigma - 1)(\sigma^2 - 1) [\Sheaf{L}]\, =\, 0& \text{if $r = 2$}.
  \end{cases}
\text{}
\eeq

where $[\Sheaf{L}]$ denotes the class of $\Sheaf{L}$ in $Pic\,E$. Regular algebras give rise to regular triples, which should come as no surprise. But, in the paper it is also shown that standard algebras give rise to regular triples [\ref{main_theorem}] and that algebras associated to regular triples are also AS-regular. 

Let $A'$ be a non-degenerate semi-standard algebra and let $\mathcal{T}(A')=(E,\sigma,\Sheaf{L})$ be the triple associated to it, with $E \hookrightarrow P^{s-1}$. Let $B=B(\mathcal{T})$ denote the twisted homogeneous coordinate ring constructed out of the triple $\mathcal{T}(A')$. Then we define another algebra $A=A(\mathcal{T}(A'))$. Set $T = \underset{n\geqslant 0}{\oplus} T_n$ be the tensor algebra over $k$ on $B_1$. The definition of a tensor algebra gives us a canonical homomorphism $T \map B$. Let $J = \underset{n\geqslant 0}{\oplus} J_n$ be its kernel, and let $I$ denote the two-sided ideal of $T$ generated by $J_s$. We then put $A(\mathcal{T})=T/I$. Then the composition of the natural homomorphisms

\beqn
A(\mathcal{T})\, =\, T/I \map T/J \map B\, =\, B(\mathcal{T})
\eeqn 

gives us a canonical homomorphism $A(\mathcal{T}) \map B(\mathcal{T})$, which is bijective in degree $1$ (since $T_1\, =\, B_1$ by definition). 

\begin{thm} \label{main_theorem}
The algebras $A$ and $A'$ are canonically isomorphic. Further, if $A'$ is standard to begin with, then $\mathcal{T}(A')$ is regular as a triple. 
\end{thm}

This theorem tells us that the AS-regular algebras of dimension $3$ are indeed the same as non-degenerate standard algebras. But non-degenerate standard algebras come hand in hand with an abstract triple and hence are more amenable to algebraic geometric tools. The story does not quite end here. Such algebras have one more desirable property. 

\begin{thm}

Let $\mathcal{T}$ be an abstract triple. Then $B=B(\mathcal{T})$ is both left and right noetherian. 

\end{thm}

Hence, AS-regular algebras of dimension $3$ are also both left and right noetherian. The proofs of these theorems are rather long and intricate. Hence, they could not be incorporated, though some ideas could be instructive. 


\newpage

\section{Towards Non-commutative Schemes via Grothendieck Categories}

\vspace{5mm}

Whether we like it or not, one has to always do some mundane things at the
beginning {\it i.e.,} fixing notations and conventions.

All categories, unless otherwise stated, will be assumed to be {\it locally small}. Not all arguments need this assumption but it just saves us some headaches. Also, since we are talking about \nc algebraic geometry, it is needless to say that rings and algebras are not assumed to be commutative. It should also be added that unadorned terms like {\it schemes} and {\it sheaves} refer to the usual (commutative) ones. 

\begin{eqnarray*}
\begin{split}
Qcoh(X) :=\; &\text{category of quasi-coherent sheaves on a scheme $X$.} \\
Coh(X)  :=\; &\text{category of coherent sheaves on $X$.} \\
Mod(A)  :=\; &\text{category of right $A$-modules, where $A$ is a $k$-algebra.} \end{split}
\end{eqnarray*}

\vspace{3mm}

Now we delve into some ``philosophy''. In mathematics we generally try to keep
things tidy and organized even if it comes for a price, say at the level of
complexity. To incorporate multifarious data we define increasingly
complicated mathematical objects, their morphisms and their compositions and
categories provide the framework where such things are practicable. Jumping
into conclusion - categories are indispensable! One might look at the set
of integers with the two operations (addition and multiplication) just as a
set or make an intellingent definition out of the whole lot of mess and call
it a {\it ring}. Now to study a ring one might try to investigate all the
questions in a myopic way looking only at the ring or map it to simple modules
and try to extract information from that. In the case of a PID the structure
theorem of the finitely generated modules over it already gives us a good idea
about the underlying structure of the ring. And now, armed with gut feeling,
we take a leap in reasoning and say that the study of a ring can be
``reduced'' to the study of the category of modules over that ring. It is not
entirely speculative as one can get a part of the ring (its centre) back by
taking the endomorphism ring of the identity functor of the module
category. As stated earlier, the category of rings is just the opposite
category of affine schemes and the category of modules is just the category of
quasi-coherent sheaves. Mathematical history reveals that whenever we have
many equivalent formulations of the same concept, we should not ignore any one
as there might be occasions where one notion generalizes, while the others do
not. So if we cannot make sense of a scheme in the direct way (it is a truism
that in the \nc world we do not have ``enough'' two-sided ideals and hence
localizations have very little chance of being pulled off), we pass on to the
level of categories. If one is not too pedantic one should not complain if the
following assertion is made - a {\it non-commutative scheme} should be an
abelian category which resembles the category of quasi-coherent sheaves on a
commutative scheme. What are such categories?? Definitely they are not unwieldy arbitrary categories; we expect them to be at least abelian. Actually what we want are Grothendieck Categories invoking the paper by Grothendieck \cite{Gro}, where he had written down the famous {\bf AB Properties} extracting geometric content out of arbitrary abelian categories. At the moment this is like putting two and two together and getting five, but hopefully this choice will be vindicated by the time we reach the end of this section. Cutting the long story short, we directly present the definition of, what we call today, a {\it Grothendieck Category}, which will be our main object of interest for some time at least. 

\vspace{3mm}
\begin{center}
{\bf Grothendieck Categories}
\end{center}

\begin{defn}
\noindent
\text{\bf Grothendieck Category}

  It is a locally small cocomplete ({\it i.e.,} closed under all coproducts)
  abelian category with a generator and satisfying, for every family of short
  exact sequences indexed by a filtered category $I$ [{\it i.e.,} $I$ is
  non-empty and, if $i,j\in Ob(I)$ then $\exists$ $k\in Ob(I)$ and arrows
  $i\map k$ and $j\map k$, and for any two arrows $i\,
  {}^{\overset{u}{\map}}_{\underset{v}{\map}}\, j$ $\exists$ $k\in Ob(I)$ and
  an arrow $w:j\map k$ such that $wu = wv$ (think of a categorical formulation
  of a directed set)]. 

\begin{eqnarray*}
0 \map A_i \map B_i \map C_i \map 0
\end{eqnarray*}

the following short sequence is also exact

\begin{eqnarray*}
0 \map \underset{\underset{i\in I}{\map}}{colim}\,A_i \map \underset{\underset{i\in I}{\map}}{colim}\,B_i \map \underset{\underset{i\in I}{\map}}{colim}\,C_i \map 0
\end{eqnarray*}

{\it i.e.,}  passing on to filtered colimits preserves exactness. This is equivalent to the sup condition of the famous {\bf AB5 Property}. 

\end{defn}

\begin{rem}

The original {\bf AB5 Property} requires the so-called {\it sup} condition, besides cocompleteness. An abelian category satisfies {\it sup} if 

for any ascending chain $\Omega$ of sub-objects of an object $M$, the supremum of $\Omega$ exists; and for any sub-object $N$ of $M$, the canonical morphism

\begin{eqnarray*}
sup\{L\cap N | L\in \Omega\} \overset{\sim}{\map} (sup\Omega) \cap N
\end{eqnarray*}

is an isomorphism. Hence, another definition of a Grothendieck Category could be a cocomplete abelian category, having a generator and satisfying the sup condition {\it i.e.,} an {\bf AB5} category with a generator. 
\end{rem}

For the convenience of the reader let me say a few things about a {\it generator} of a category. An object $G$ of a category $\mathcal{C}$ is called a {\it generator} if, given a pair of morphisms $f,g:A \map B$ in $\mathcal{C}$ with $f\neq g$, there exists an $h:G\map A$ with $fh\neq gh$ (more crisply, $Hom(G,{}_{-}):\mathcal{C} \map Set$ is a faithful functor). A family of objects $\{G_i\}_{i \in I}$ is called a {\it generating set} if, given a pair of morphisms $f,g:A{}_{\map}^{\map} B$ with $f \neq g$, there exists an $h_i:G_i \map A$ for some $i \in I$ with $fh_i \neq gh_i$. Strictly speaking, this is a misnomer. In a cocomplete category, a family of objects $\{G_i\}_{i\in I}$ forms a generating set if and only if the coproduct of the family forms a generator. 

\begin{rem}
Let $\mathcal{C}$ be a cocomplete abelian category. Then an object $G$ is a generator if and only if, for any object $A \in Ob(\mathcal{C})$ there exists an epimorphism, 

\begin{eqnarray*}
G^{\oplus I} \map A
\end{eqnarray*} 

for some indexing set $I$. 

\end{rem}

\begin{ex}
\textbf{(Grothendieck Categories)}\newline

1. $Mod(R)$, where $R$ is an associative ring with unity. 

2. The category of sheaves of $R$-modules on an arbitrary topological space.

3. In the same vein, the category of abelian pre-sheaves on a {\it site} $\mathcal{T}$. Actually this is just $Funct(\mathcal{T}^{op},Ab)$.

4. $Qcoh(X)$, where $X$ is a quasi-compact and quasi-separated scheme. (a morphism of schemes $f:X\map Y$ is called quasi-compact if, for any open quasi-compact $U\subseteq Y$, $f^{-1}(U)$ is quasi-compact in $X$ and it becomes quasi-separated if the canonical morphism $\delta_f:X\map X\times_Y X$ is quasi-compact. A scheme $X$ is called quasi-compact (resp. quasi-separated) if the canonical unique morphism $X\map Spec(\ZZ)$, $Spec(\ZZ)$ being the final object, is quasi-compact (resp. quasi-separated)).    
\end{ex}

Grothendieck categories have some remarkable properties which make them amenable to homological arguments.

1. Grothendieck categories are complete {\it i.e.,} they are closed under products.  

2. In a Grothendieck category every object has an injective envelope, in particular there are enough injectives.

\vspace{2mm}

The deepest result about Grothendieck categories is given by the following theorem. 

\begin{thm} \text{[Gabriel Popescu]} \label{Gabriel_Popescu}
Let $\mathcal{C}$ be a Grothendieck category and let $\mathcal{G}$ be a generator of $\mathcal{C}$. Put $S = End(\mathcal{G})$. Then the functor 

\begin{eqnarray*}
Hom(\mathcal{G},{}_{-}):\mathcal{C}\map Mod(S^{op})
\end{eqnarray*}

is fully faithful (and has an exact left adjoint). 
\end{thm}

We wrap things up in this episode by showing that Grothendieck categories are, in some sense, ``big'' amongst abelian categories. Let $\mathcal{C}$ be a small abelian category. Then define the category $Ind\mathcal{C}$ as:

\begin{eqnarray*} \label{IndC}
\begin{split}
Ob(Ind\mathcal{C}) = &\text{Functors from all filtered categories to $\mathcal{C}$.} \\
                     &\text{For convenience we denote them $\{T_j\}_{j\in J}$, where $J$ is a filtered category.}
\end{split}
\end{eqnarray*}

\noindent
$Hom_{Ind\mathcal{C}}(\{T_i\}_{i\in I},\{T'_j\}_{j\in J})\,=\, Lim_{j\in J}\,colim_{i\in I}\,Hom_{\mathcal{C}}(T_i,T'_j)$

\vspace{1mm}

By means of a dual construction one arrives at, what is called, a pro-category $Pro\mathcal{C}$. The category $Ind\mathcal{C}$ has a more tangible realisation in the form of a naturally equivalent category, denoted $Lex(\mathcal{C})$, where \\
$Ob(Lex(\mathcal{C}))\,=\,\{\text{left exact additive functors $\mathcal{C}^{op} \map Ab$}\}$. By {\it Yoneda's theorem} the functor 

\begin{eqnarray*}
\begin{split}
&\mathcal{C} \map Lex(\mathcal{C}) \\
& A \longmapsto Hom_{\mathcal{C}}({}_{-},A)
\end{split}
\end{eqnarray*}

is fully faithful and exact. This gives an embedding of the abelian category $\mathcal{C}$ into $Lex(\mathcal{C})$. It is known that if $\mathcal{C}$ is a small abelian category and $\mathcal{D}$ is a Grothendieck category, then $Funct(\mathcal{C}^{op},\mathcal{D})$ is also a Grothendieck category. Hence, in particular, $Lex(\mathcal{C})$ is a Grothendieck category. So we obtain an embedding of a small abelian category into a Grothendieck category. Now we mention a theorem which will come in handy later on. 

\begin{thm} \text{[P. Gabriel]} \label{noeVsLocNoe}
Let $\mathcal{C}$ be a noetherian abelian category. Then the categories $\mathcal{C}$ and $Ind{\mathcal{C}}$ determine each other up to a natural equivalence. In fact, $\mathcal{C}$ is the full sub-category formed by all the noetherian objects in $Ind{\mathcal{C}}$.
\end{thm}

If $X$ is a noetherian scheme, then $Ind(Coh(X))\simeq Qcoh(X)$. Let $Noeth$ denote the operation of taking noetherian objects. Then, $Noeth(Qcoh(X))\simeq Coh(X)$. So for a noetherian scheme $X$, we might as well deal with the more tractable category $Coh(X)$.

\begin{center}
{\bf Justification for bringing in Grothendieck Categories}
\end{center}

\vspace{2mm}

We begin by directly quoting {\it Manin} \cite{Man2} - ``...Grothendieck taught us, to do geometry you really don't need a space, all you need is a category of sheaves on this would-be space.'' This idea gets a boost from the following reconstruction theorem. 

\begin{thm}\text{[Gabriel Rosenberg \cite{Ros2}]}  \label{reconstruction}
Any scheme can be reconstructed uniquely up to isomorphism from the category of quasi-coherent sheaves on it.
\end{thm}

\noindent
{\bf An easy reconstruction:}

Let us look at a very simple case where this theorem is applicable. When it is known in advance that the scheme to be reconstructed is an affine one, we can just take the centre of the category, which is the endomorphism ring of the identity functor of the category. More precisely, let $X = Spec\,A$ be an affine scheme and let $\mathcal{A}$ be the category of quasi-coherent sheaves on $X$, which is the same as $Mod(A)$. Then the centre of $\mathcal{A}$, denoted $End(Id_\mathcal{A})$, is canonically isomorphic to $A$. [Centre of an abelian category is manifestly commutative and, in general, it gives us only the centre of the ring, that is $\mathcal{Z}(A)$. But here we are talking about honest schemes and hence, $\mathcal{Z}(A)=A$]. 

\begin{eqnarray*}
\begin{split}
\psi: A &\map End(Id_\mathcal{A}) \\
      a &\longmapsto \psi_a\; \text{such that $\psi_a:Id_\mathcal{A}(M)\map Id_\mathcal{A}(M)$ is just mult. by $a$ $\forall$ $M\in Mod(A)$.}
\end{split}
\end{eqnarray*} 

It is easy to see that the map is injective and a ring homomorphism. Now any $\theta \in End(Id_\mathcal{A})$ is a collection of endomorphisms $\{\theta_M\}_{M\in Mod(A)}:M\map M$, such that for all $M,N\in Mod(A)$ and for all $\phi\in Hom(M,N)$,

\[
\begin{CD}
M @>{\theta_M}>> M \\
@V{\phi}VV @VV{\phi}V \\
N @>>{\theta_N}> N
\end{CD}
\]

the above diagram commutes. In particular, when $M=A$, using $Hom_A(A,M)\cong M$ via $\phi\mapsto \phi(1)$, we see that $\theta_N\circ\phi(1) = \phi\circ\theta_A(1) = \theta_A(1).\phi(1)$. So $\theta_N =$ mult. by $\theta_A(1)\in A$ independent of $N$.

\begin{flushright}
$\blacksquare$
\end{flushright}

We can also get a derived analogue of the above result, which is, however, considerably weaker. Also it is claimed to be an easy consequence of the above theorem in \cite{BO1}.
 
\begin{thm} \text{[Bondal Orlov \cite{BO1}]}
Let X be a smooth irreducible projective variety with ample canonical or anti-canonical sheaf. If $\mathcal{D} = D^b Coh(X)$ is equivalent as a graded category to $D^b Coh(X')$ for some other smooth algebraic variety $X'$, then $X$ is isomorphic to $X'$. 
\end{thm}

Finally, consider a pre-additive category with a single object, say $\ast$. Then being a pre-additive category $Hom(\ast,\ast)$ is endowed with an abelian group structure. If we define a product on it by composition, then it is easy to verify that the two operations satisfy the ring axioms. So $Hom(\ast,\ast)$ or simply $End(\ast)$ is a ring and that is all we need to know about the pre-additive category. Extrapolating this line of thought, we say that pre-additive categories generalize the concept of rings and since schemes are concocted from commutative rings, it is reasonable to believe that by some suitable constructions on something like a pre-additive category with some geometric properties ($=$ a Grothendieck category) we can find ``\nc schemes''. 
\vspace{3mm}
\begin{center}
{\bf A small discussion on construction of quotient categories}
\end{center}

\vspace{2mm}

Recall that we call a full sub-category $\mathcal{C}$ of an abelian category $\mathcal{A}$ {\it thick} if the following condition is satisfied:

\vspace{1mm}

\hspace{4mm} for all short exact sequences in $\mathcal{A}$ of the form $0\map M'\map M\map M''\map 0$, we have

\begin{eqnarray*}
M\in \mathcal{C} \Longleftrightarrow \text{ both } M',M''\in \mathcal{C}
\end{eqnarray*}

\vspace{1mm}

It becomes a \underline{\it Serre sub-category} if it is {\it coreflective}
({\it i.e.,} the canonical inclusion functor $\mathcal{C}\map\mathcal{A}$ has
a right adjoint), besides being {\it thick}. This is essentially a ``local
terminology''. 

\begin{rem}
What we call thick sub-categories out here are referred to as Serre sub-categories by some authors. Further, note that in this definition, Serre sub-categories are themselves cocomplete (the image of the coproduct in the bigger category under the adjoint functor of the inclusion belongs to the Serre sub-category and is easily seen to be a coproduct as well). 
\end{rem}  

Now we construct the quotient of $\mathcal{A}$ by a thick sub-category $\mathcal{C}$, denoted $\mathcal{A/C}$, as follows: 

\begin{eqnarray*}
\begin{split}
Ob(\mathcal{A/C}) &= \text{objects of $\mathcal{A}$.} \\
Hom_{\mathcal{A/C}}(M,N) &= \underset{\underset{\underset{s.t.\, M/M',N'\in\mathcal{C}}{\forall\,M',N'\,sub-obj.}}{\longrightarrow}}{Lim}Hom_{\mathcal{A}}(M',N/N').
\end{split}
\end{eqnarray*}

Of course, one needs to check that as $M'$ and $N'$ run through all sub-objects of $M$ and $N$ respectively, such that $N'$ and $M/M'$ are in $Ob(\mathcal{C})$, the abelian groups $Hom_{\mathcal{C}}(M',N/N')$ form a directed system. 

The essence of this quotient construction is that the objects of $\mathcal{C}$ become isomorphic to zero. 

\begin{ex}
Let $\mathcal{A} = Mod(\ZZ)$ and $\mathcal{C} = Torsion\;groups$.
Then one can show that $\mathcal{A}/\mathcal{C} \simeq Mod(\QQ)$.

Let us define a functor from $\mathcal{A}/\mathcal{C}$ to $Mod(\QQ)$
by tensoring with $\QQ$. We simplify the $Hom$ sets of
$\mathcal{A}/\mathcal{C}$. Using the structure theorem, write every
abelian group as a direct sum of its torsion part and torsion-free
part. If one of the variables is torsion, it can be shown that in
the limit $Hom$ becomes $0$. So we may assume that both variables
are torsion-free and for simplicity let us consider both of them to
be $\ZZ$. Then,

\beqn
\begin{split}
Hom_{\mathcal{A}/\mathcal{C}}(\ZZ,\ZZ)&=\underset{n}{{\underset{\map}{{Hom}_{\mathcal{A}}}}} (n\ZZ,\ZZ) \\
                             &=\underset{n}{\bigcup}\;\frac{1}{n}\ZZ \\
                             &=\QQ = Hom(\QQ,\QQ)
\end{split}
\eeqn

This says that the functor is full, and an easy verification shows that it is faithful and essentially surjective.

\end{ex}

Let $\pi$ denote the canonical quotient functor $\mathcal{A}\map\mathcal{A/C}$. If $\pi$ has a right adjoint $S:\mathcal{A/C}\map\mathcal{A}$  called the {\it section functor}, then the thick sub-category $\mathcal{C}$ is called {\it localizing}. It is known that in a Grothendieck category, the notions of a Serre sub-category and a localizing sub-category are the same. Further, if $\grot$ is a Grothendieck category and $\mathcal{C}$ is a localizing sub-category, then $\grot/\mathcal{C}$ is again a Grothendieck category. 

\vspace{2mm}
\begin{center}
{\bf Spectrum of Grothendieck Categories}
\end{center}

The passage between commutative rings and affine spaces is so smooth and transparent that one gets easily tempted to directly generalize such ideas to the \nc domain. For \nc rings which are finite modules over their centres, there exists a satisfactory theory. Let $A$ be a \nc ring with centre, $\mathcal{Z}(A)=R$. Then the pair $(Spec\,R,\tilde{A})$ is a nice geometrical object, where $\tilde{A}$ is the sheaf of algebras associated to the $R$-algebra $A$, satisfying $Mod(A)\simeq$ category of sheaves of $\tilde{A}$-modules. But many interesting algebras encountered in physics fail to be centre-finite modules, for example, {\it Weyl algebras}. To associate a ``space'' to a \nc ring we work with not just the ring, but the category of modules over that ring. Don't forget our paradigm - a ``space'' is just the category of sheaves on it. But such a category is a Grothendieck category. So we are ready to enter the murky world of Grothendieck categories. The contents under this sub-title are based entirely upon chapter 3 of {\it A. Rosenberg's} book \cite{Ros1} and for details interested readers are requested to look it up (he talks about arbitrary abelian categories but we focus on Grothendieck categories). The aim is to concoct something resembling a scheme from a Grothendieck category. We shall discuss the different attributes of a scheme very tersely one by one. 

\vspace{2mm}
\noindent
\underline{\bf The underlying set ($Spec(\grot)$)}

\begin{defn}
{\bf A pre-order on Grothendieck Categories ($\succ$):} 

Fix a Grothendieck Category $\grot$. For any two objects $X$ and $Y$ of $\grot$ we shall write $X\succ Y$ if and only if $Y$ is a subquotient of a finite coproduct of copies of $X$ {\it i.e.,} there exists $U$ such that the following holds

\begin{eqnarray*}
\underset{finite}{\oplus}X\hookleftarrow U \twoheadrightarrow Y
\end{eqnarray*}

\end{defn}

A pre-order is just reflexive and transitive. We introduce an equivalence and go modulo that equivalence to bring in symmetry. Call $X \sim Y$ if and only if $X\succ Y$ and $Y\succ X$. Then set $|\grot| = ((Ob(\grot)/{\sim})\, ,\,\succ)$. This is an ordered set. 

\begin{defn}
{\bf Spectrum of $\grot$:} 

Let $S(\grot) = \{M\in Ob(\grot)\,|\, M\neq 0, \forall\, 0\neq N\hookrightarrow M,\,N\succ M\}$. The spectrum of $\grot$, denoted $Spec(\grot)$, is the ordered set of equivalence classes (with respect to $\succ$) of elements of $S(\grot)$. 
\end{defn}

\begin{rem}
It is obvious that {\it simple} objects (objects which have no non-zero proper sub-objects) belong to $S(\grot)$. It can also be shown that two simple objects are equivalent (with respect to $\succ$) if and only if they are isomorphic. So in $Spec(\grot)$ isomorphic simple objects are clubbed together. 
\end{rem}

\begin{prop}
Let $Q:\grot \map \mathcal{B}$ be an exact localization functor between Grothendieck categories. Then, for any $P\in Spec(\grot)$, either $Q(P)=0$ or $Q(P)\in Spec(\mathcal{B})$. 
\end{prop} 

This proposition says that exact localizations almost respect spectrums ($0$ does not belong to the spectrum). 

\vspace{2mm}
\noindent
\underline{\bf Realisation of $Spec(\grot)$ in terms of Serre Sub-categories}

\vspace{2mm}

While discussing quotient construction in categories we have already introduced the notion of a {\it Serre Sub-category}. 

\vspace{2mm}
\noindent
{\bf A useful notation $\langle M\rangle$:}

For any $M\in Ob(\grot)$ put

\begin{eqnarray*}
 \langle M\rangle\, :=\,\text{full sub-category of $\grot$ generated by objects $N$, such that $N\nsucc M$.}
\end{eqnarray*}
 
\begin{lem}
For any $M,M'\in Ob(\grot)$, $M\succ M'$ if and only if $\langle M\rangle \supseteq \langle M'\rangle$. Thus, $M$ is equivalent to $M'$ with respect to $\succ$ if and only if $\langle M'\rangle = \langle M\rangle$. 
\end{lem}

Let us continue to denote by $\langle M\rangle$ the equivalence class of $\langle M\rangle$. Then we have the following map between posets. 

\begin{eqnarray*}
\begin{split}
(|\grot|,\succ) &\map (\{\langle M\rangle | M\in Ob(\grot)\},\supseteq) \\
 M              &\longmapsto \langle M\rangle
\end{split}
\end{eqnarray*}

This is clearly one-to-one and hence, we have a realisation of $|\grot |$ in terms of certain sub-categories of $\grot$. 

\begin{prop} \label{realisation1}
If an object $P$ of the category $\grot$ belongs to $Spec(\grot)$, then $\langle P\rangle$ is a Serre sub-category of $\grot$.
\end{prop}

\vspace{2mm}
\noindent
{\bf Sketch of Proof:}

Let $P\in Spec(\grot)$. We need to show that for any exact sequence $0\map M'\map M\map M''\map 0$, $M\in \langle P\rangle$ if and only if both $M'\in \langle P\rangle$ and $M''\in \langle P\rangle$. Consider the following diagram,

\begin{eqnarray*}
\begin{split}
0 & \\
\uparrow & \\
\underset{finite}{\oplus}M'' & \\
\uparrow & \overset{i}{\nwarrow} \\
\underset{finite}{\oplus}M & \overset{j}{\hookleftarrow} K \twoheadrightarrow P \\
\uparrow & \underset{k}{\swarrow} \\
\underset{finite}{\oplus}M' & \\
\uparrow & \\
0
\end{split}
\end{eqnarray*}

where $i,j$ and $k$ are supposed to be injections. Then existences of $i$ and $k$ easily imply the existence of $j$ (just compose $k$ with the map $\underset{finite}{\oplus}M' \map \underset{finite}{\oplus}M$). Now suppose $j$ exists. Then, $\underset{finite}{\oplus}M' \hookleftarrow K\cap (\underset{finite}{\oplus}M') \map P$ $\Longrightarrow$ $M'\succ L\hookrightarrow P$ and $\underset{finite}{\oplus}M''\hookleftarrow K' \map P$ $\Longrightarrow$ $M''\succ L'\hookrightarrow P$, where $K'=K/(K\cap (\underset{finite}{\oplus}M'))$. But as $P\in Spec(\grot)$, $L,L'$ sub-objects of $P$ $\Rightarrow$ $L,L'\succ P$. Hence, both $M',M''\succ P$. We have actually proved, $M\notin \langle P\rangle \Longleftrightarrow \text{ both }M',M'' \notin \langle P\rangle$. 

For coreflectiveness, we need a right adjoint for $\imath:\langle P\rangle\map\grot$. For any $M\in Ob(\grot)$, let $\langle P\rangle (M)$ denote the set of all sub-objects of $M$ which belong to $\langle P\rangle$. Since $\grot$ is a Grothendieck category it has the (sup) property. So $\bar{M}:=\,sup(\langle P\rangle (M))$ exists in $\grot$. But, as $\langle P\rangle$ contains all sub-objects and quotient-objects of its objects, $\bar{M}\in Ob(\langle P\rangle)$. Define the adjoint of $\imath$, denoted $\imath^{!}$, by 

\begin{eqnarray*}
\begin{split}
\imath^{!}:\grot &\map\langle P\rangle \\
             M   &\longmapsto \bar{M}
\end{split}
\end{eqnarray*}

\begin{flushright}
$\blacksquare$
\end{flushright}

It requires some more work to prove the following proposition.

\begin{prop} \label{realisation2}
For any object $M$ of $\grot$, such that $\langle M\rangle$ is a Serre sub-category, there is an object $P\in Spec(\grot)$ which is equivalent to $M$, {\it i.e.,} $\langle M\rangle = \langle P\rangle$. 
\end{prop}

Putting two [\ref{realisation1}] and two [\ref{realisation2}] together we get a realisation of $Spec(\grot)$.

\vspace{1mm}
\begin{center}
$Spec(\grot)$ $=$ \{$\langle M\rangle$ $|$ $M\in Ob(\grot)$ and $\langle M\rangle$ is a Serre sub-category\}
\end{center}

\vspace{3mm}
\noindent
\underline{\bf Categorical incarnation of Local rings of $Spec(\grot)$}
 
\vspace{2mm}

As we have just seen, to every object $M\in Spec(\grot)$ one can associate a localizing sub-category $\langle M\rangle$ and hence an exact localization $Q_M:\grot\map\grot/\langle M\rangle$. 

A non-zero object $M$ of $\grot$ is called {\it quasi-final} if $N\succ M$ for any non-zero object of $\grot$ or equivalently $\langle M\rangle$ is the zero sub-category. 

A quasi-final object, if it exists, clearly belongs to the spectrum and any two such objects are equivalent.

\begin{defn} {\bf Local abelian category}
An abelian category $\mathcal{A}$ is called local if it has a quasi-final object. 
\end{defn}

Now we have the following proposition to justify the christening.

\begin{prop} \label{local}
Let $\mathcal{A}$ be a local abelian category. Then the center of $\mathcal{A}$ is a commutative local ring.
\end{prop}

Though we have introduced the notion of locality for arbitrary abelian categories, we go back to our favourite Grothendieck categories. In a usual (commutative) scheme the stalk of the structure sheaf at a point is a local ring. So far we have obtained a set, which is $Spec(\grot)$, and if we have done something sane then for every point $M\in Spec(\grot)$ we should expect to get a ``local object''. Coupled with [\ref{local}] the following proposition gives us precisely what we expect.

\begin{prop}
For any object $M\in Spec(\grot)$, the quotient category $\grot/\langle M\rangle$ is local.
\end{prop}

\noindent
{\bf Sketch of proof:}

Let $Q:\grot\map\grot/\langle M\rangle$ be the exact localization functor. We shall show that $Q(M)$ is quasi-final and hence, $\grot/\langle M\rangle$ is local. Take any non-zero object $X\in \grot/\langle M\rangle$. $X$ can also be thought of as an object in $\grot$ [refer to quotient construction] and since it is a non-zero object in the quotient category, $X\notin Ob(\langle M\rangle)$. This implies that $X\succ M$. Since $Q$, being an exact functor, respects $\succ$, we have $Q(X)=X\succ Q(M)$. Since $X$ was chosen arbitrarily, $Q(M)$ is a quasi-final object. 

\begin{flushright}
$\blacksquare$
\end{flushright}

\vspace{3mm}
\noindent
\underline{\bf Categorical incarnation of the residue field at a point}
\vspace{2mm}

As is wont, we need to beat around the bush a little bit. Some more technical jargon have to show up.

\begin{defn} {\bf Support of an object:}

The support of $M\in Ob(\grot)$, denoted $Supp(M)$ is the set of all $\langle P\rangle\in Spec(\grot)$ such that $M\succ P$. Or, in other words,

\begin{center}
$Supp(M) = \{\langle P\rangle\in Spec(\grot) | Q_{\langle P\rangle}M\neq 0\}$
\end{center}

\end{defn}

For any subset $W\subseteq Spec(\grot)$, let $A(W)$ be the full sub-category generated by all objects $M$ of the category $\grot$ such that $Supp(M)\subseteq W$. 

\begin{lem} \label{char_of_A(W)}

$A(W) = \underset{P\in W^{\perp}}{\cap} \langle P\rangle$, where $W^{\perp}:= Spec(\grot) - W$.

\end{lem} 

\noindent
{\bf Sketch of proof:}

Let $M\in Ob(\underset{P\in W^{\perp}}{\cap}\langle P\rangle)$ $\Longleftrightarrow$ $Supp(M)\cap Ob(\underset{P\in W^{\perp}}{\cap} \langle P\rangle) = \emptyset$ $\Longleftrightarrow$ $Supp(M)\subseteq W$ {\it i.e.,} $M\in Ob(A(W))$. 

\begin{flushright}
$\blacksquare$
\end{flushright}

Now if one is ready to believe that the intersection of any set of Serre sub-categories is once again a Serre sub-category, we may conclude from the above lemma that $A(W)$ is a Serre sub-category. For any point $\langle P\rangle\in Spec(\grot)$ consider the category 

\begin{eqnarray*}
\bar{\mathcal{K}}(P) := A(Supp(P))/\langle P\rangle
\end{eqnarray*}

Now we denote by $\mathcal{K}(P)$ the full sub-category of $\bar{\mathcal{K}}(P)$ generated by all objects $M$ of $\bar{\mathcal{K}}(P)$ which are supremums of their sub-objects V, which satisfy $\langle V\rangle =\langle P\rangle$. 

We call the category $\mathcal{K}(P)$ the {\it residue category} of $\langle P\rangle$. One can check that it is local and its spectrum is singleton. 
 
\vspace{3mm}
\noindent
\underline{\bf The Zariski topology on $Spec(\grot)$}
\vspace{2mm}

To define the analogue of Zariski topology we need some ingredients first. 

\begin{defn} {\bf Topologizing sub-category:}

A full sub-category $\mathcal{C}$ of $\grot$ is called topologizing if for any exact sequence in $\grot$ 

\begin{eqnarray*}
0\map M'\map M\map M''\map 0
\end{eqnarray*}

$M\in \mathcal{C}$ $\Longrightarrow$ both $M',M''\in \mathcal{C}$. Besides, it should also be closed under finite coproducts.
\end{defn}

Clearly any thick sub-category is topologizing, and so is any Serre sub-category (Serre $\Longrightarrow$ thick).  

\begin{defn} {\bf Gabriel product of topologizing sub-categories:}

For any two topologizing sub-categories $\mathcal{S},\mathcal{T}$ of $\grot$ define the Gabriel product $\mathcal{S}\bullet\mathcal{T}$ as the full sub-category of $\grot$ generated by all objects $M$ of $\grot$ such that there exists an exact sequence in $\grot$

\begin{eqnarray*}
0\map M'\map M\map M''\map 0
\end{eqnarray*}

with $M'\in Ob(\mathcal{T})$ and $M''\in Ob(\mathcal{S})$. 
\end{defn}

One can check that the Gabriel product of two topologizing sub-categories is also a topologizing one. Further,

\begin{center}
$\mathcal{S}\bullet (\mathcal{T}\bullet\mathcal{U}) = (\mathcal{S}\bullet\mathcal{T})\bullet\mathcal{U}$ and $0\bullet\mathcal{S} = \mathcal{S} = \mathcal{S}\bullet 0$
\end{center}

It is also evident from the definitions that a topologizing sub-category $\mathcal{T}$ is thick if and only if $\mathcal{T}\bullet\mathcal{T} = \mathcal{T}$.

\begin{defn} {\bf Closed and left closed sub-categories:}

A topologizing sub-category is closed (resp. left closed) if is coreflective (resp. reflective) as well. [coreflective (resp. reflective) means that the canonical inclusion functor has a right (resp. left) adjoint].
\end{defn}

It needs some formal arguments to prove that the Gabriel product of two closed (resp. left closed) sub-categories is also closed (resp. left closed). Left closedness is also preserved under arbitrary intersections. Now we define for any topologizing sub-category (in particular, for any left closed sub-category)

\vspace{1mm}
\noindent
$V(\mathcal{T})$ $:=$ $\{\langle P\rangle\in Spec(\grot) | P\in Ob(\mathcal{T})\}$

This notation is not without any foresight. Now we are ready to define the Zariski Topology, denoted $\mathcal{ZT}$. 

\begin{defn} {\bf Zariski Topology ($\mathcal{ZT}$):}

It is the topology generated by the closed sets of the form $V(\mathcal{T})$, where $\mathcal{T}$ runs through the set of all left closed sub-categories of $\grot$. 
\end{defn}

The fact that the sets of the form $V(\mathcal{T})$ play the role of the closed subsets of $\mathcal{ZT}$ is corroborated by the following proposition. 

\begin{prop}
The set $\mathcal{ZT}$ is closed under finite unions and arbitrary intersections.
\end{prop}

\noindent
{\bf Sketch of proof:}

(a) Union: $V(\mathcal{S})\cup V(\mathcal{T}) = V(\mathcal{S}\bullet\mathcal{T})$. 

Clearly $\mathcal{S}\subseteq\mathcal{S}\bullet\mathcal{T}\supseteq\mathcal{T}$. which implies the inclusion $V(\mathcal{S})\cup V(\mathcal{T})\subseteq V(\mathcal{S}\bullet\mathcal{T})$. 

For the other inclusion let us take any $\langle P\rangle\in V(\mathcal{S}\bullet\mathcal{T})$ {\it i.e.,} $P\in Spec(\grot)\cup Ob(\mathcal{S}\bullet\mathcal{T})$. The latter means that there exists an exact sequence 

\begin{center}
$0\map P'\map P\map P''\map 0$
\end{center}

in which $P'\in Ob(\mathcal{T})$ and $P''\in Ob(\mathcal{S})$. 

If $P'\neq 0$, then $P'\succ P$ $\Longrightarrow$ $P\in Ob(\mathcal{T})$.

If $P'=0$, then $P\simeq P''$ $\Longrightarrow$ $P\in Ob(\mathcal{S})$. 

\noindent
Therefore, we have

\begin{center}
$Spec(\grot)\cup Ob(\mathcal{S}\bullet\mathcal{T})\subseteq (Spec(\grot)\cup Ob(\mathcal{S}))\cup(Spec(\grot)\cup Ob(\mathcal{T}))$
\end{center}

which implies the desired inclusion.

\vspace{2mm}
\noindent
(b) Intersection: 

Let $\Omega$ be any set of left closed sub-categories of $\grot$. Clearly, $\underset{\mathcal{T}\in\Omega}{\cap} V(\mathcal{T}) = V(\underset{\mathcal{T}\in\Omega}{\cap} \mathcal{T})$. But left closed sub-categories, as stated earlier, are closed under arbitrary intersections. So $\underset{\mathcal{T}\in\Omega}{\cap}\mathcal{T}$ is also left closed $\Longrightarrow$ $V(\underset{\mathcal{T}\in\Omega}{\cap}\mathcal{T})\in\mathcal{ZT}$.

\begin{flushright}
$\blacksquare$
\end{flushright}

\begin{rem}

It has been shown in the book \cite{Ros1} that $(Spec(\grot),\mathcal{ZT})$ is quasi-compact as a topological space if the generator of the Grothendieck Category $\grot$ is of finite type ({\it an object $M$ of $\grot$ is said to be of finite type if, for any directed set $\Omega$ of sub-objects of $M$ such that $sup\,\Omega = M$, there already exists a sub-object $M'\in\Omega$ such that $M'=M$}). Further, for an associative unital ring $R$, $(Spec(Mod(R)),\mathcal{ZT})$ is quasi-compact and has a basis of quasi-compact open subsets. It should also be bourne in mind that the Zariski topology is not satisfactory in all circumstances. There are some other canonical topologies, denoted $\tau^\ast$ and $\tau_\ast$, for which interested readers may refer to \cite{Ros1} and \cite{Ros2}. 

\end{rem}

\vspace{3mm}
\noindent
\underline{\bf Sheaves on $Spec(\grot)$}
\vspace{3mm}

We already know that the closed sets of $\mathcal{ZT}$ are of the form $V(\mathcal{T})$. Taking the complements in $Spec(\grot)$ we get the open sets, which we denote by $Open(\mathcal{ZT})$ (this can easily be thought of as a category with inclusions being the only morphisms). To any $U\in Open(\mathcal{ZT})$ we assign the Serre sub-category $\langle U\rangle := \underset{P\in U}{\cap}\langle P\rangle$. According to lemma [\ref{char_of_A(W)}] $\langle U\rangle = A(U^\perp)$, where $U^\perp$ is just the complementary closed subset of $U$ in $\mathcal{ZT}$. 

The ``structure sheaf'' is obtained as a sheaf of commutative rings via the functor $\sheaf_\grot$.

\begin{eqnarray*}
\begin{split}
\sheaf_\grot:Open(\mathcal{ZT})^{op} &\map (Comm.\;rings) \\
                 U         &\longmapsto \mathcal{Z}(\grot/\langle U\rangle)
\end{split}
\end{eqnarray*}

where $\mathcal{Z}(\grot/\langle U\rangle)$ stands for the centre of the category $\grot/\langle U\rangle$. 

\vspace{2mm}
\noindent
To any $M\in Ob(\grot)$ we associate a functor $\tilde{M}$. Let $Q_{\grot/\langle U\rangle}:\grot\map\grot/\langle U\rangle$ and put $M_U = Q_{\grot/\langle U\rangle}(M)$. 

\begin{eqnarray*}
\begin{split}
\tilde{M}:Open(\mathcal{ZT})^{op} \map & Ab \\
              U              \longmapsto & M_U \\
         & \text{(in an ab. cat. every object is an abelian group object.)}
\end{split}
\end{eqnarray*}

The $\tilde{M}$, so defined, is called the {\it sheaf associated to the object $M$}. Of course, one has to check that both $\sheaf_\grot$ and $\tilde{M}$ satisfy the sheaf conditions. Some more work actually shows that $\tilde{M}$ defines a sheaf of $\sheaf_\grot$-modules ($\theta:Id_{\grot/\langle U\rangle}\map Id_{\grot/\langle U\rangle}$ then $\theta$ is easily seen to act on $M_U$). Further, one has $Mod(\sheaf_\grot(U))\simeq \grot/\langle U\rangle$.

In the book \cite{Ros1} {\it Rosenberg} also talks about relative spectra associated to functors $F:\mathcal{A}\map\mathcal{B}$, but due to lack of time and space we skip those details. 

\vspace{3mm}
\noindent
{\bf A glaring omission:}

{\it A. Rosenberg} has argued in his book \cite{Ros1} that his spectrum is actually smaller than the other prevalent ones and is also easier to compute. As an illustration he has computed the spectra of so-called {\it hyperbolic rings} {\it i.e.,} rings of the form $R\langle X,Y\rangle$ subject to the relations $Xr=\theta(r)X$, $Yr=\theta^{-1}(r)Y$, $XY=u$ and $YX=\theta^{-1}(u)$, where $\theta \in Aut(R)$,$r\in R$ and $u\in \mathcal{Z}(R)$. Such a spectrum, when specialized to a module category over an associative unital ring $R$, coincides with the so-called {\it left spectrum} of $R$. 

\vspace{2mm}
\noindent
{\bf Conclusion:}

We have pushed a heavy ball up to the top of a hill and now we can just let it roll down via the reconstruction theorem [\ref{reconstruction}] to get back to the square position. To be precise, suppose $\grot = Qcoh(X)$, where $X$ is a noetherian (actually quasi-compact and quasi-separated is enough) scheme. Then it has been proved in \cite{Ros2} that the ringed space $(Spec(\grot),\sheaf_\grot)$ is isomorphic to $(X,\sheaf_X)$. Actually the reconstruction works for arbitrary schemes. We have deliberately worked over Grothendieck categories to ``conform to the title of the section''. Unfortunately, the category of quasi-coherent sheaves may not be a Grothendieck category for arbitrary schemes, but for noetherian ones, it is. For a noetherian scheme $X$, we are allowed to concentrate on $\grot = Coh(X)$ by [\ref{noeVsLocNoe}] and the comment following that. Then the reconstruction maps are (in less than a nutshell)

\begin{eqnarray*}
\begin{split}
\phi:X \map & Spec(\grot) \\
     x \longmapsto &\text{$\langle P_x\rangle$ the sky-scraper sheaf,} \\
            &\text{whose only non-zero stalk in the residue field at $x$.}
\end{split}
\end{eqnarray*}

\begin{eqnarray*}
\begin{split}
\psi:Spec(\grot) \map & X \\
\langle M\rangle \longmapsto & \text{the generic point of the support of $M$.} \\
       & \text{$\langle M\rangle\in Spec(\grot)$ $\Longrightarrow$ $Supp(M)$ is irreducible and closed in $X$.}
\end{split}
\end{eqnarray*}

One has to check that $\psi\circ\phi = Id_X$ and $\phi\circ\psi = Id_{Spec(\grot)}$ and also that they are indeed morphisms of ringed spaces. The definition of $Spec(\grot)$ is such that sheaves with reducible supports cannot belong to it. This is one of the interesting points of the proof.   
 
This excerpt is in the spirit of \nc algebraic geometry subsuming representation theory. The classical problem of determining all irreducible representations of an algebra $R$ (typically enveloping algebras) can be reduced to the study of the spectrum of the category, $Mod(R)$. We have provided a general recipe to construct ringed spaces from Grothendieck categories and $Mod(R)$ is an archetype of such categories. The existence of a generator in a Grothendieck category is quite a strong hypothesis as we know virtually everything about the category from the module category of the endomorphism ring of the generator ({\it Gabriel-Popescu} [\ref{Gabriel_Popescu}]). This \nc spectral theory reveals important classes of irreducible representations which are, otherwise, out of reach of the classical theory ({\it Verma and Harish-Chandra modules}). 


\newpage

\section{Non-commutative Projective Geometry}

\vspace{5mm}

     Although so far we have done everything in a very general set-up, we would soon like to concentrate on, what is called, {\it non-commutative projective geometry}. But before that let us go through one nice result in the affine case. Let $X$ be an affine scheme and put $A = \Gamma(X,\sheaf_X)$. Then it is classical that $Qcoh(X)$ is equivalent to $Mod(A)$. This immediately begs the question of which Grothendieck categories can be written as $Mod(A)$ for some possibly \nc ring $A$. The answer is given by the theorem below.

\begin{thm} \text{\cite{Ste}}
Let $\mathcal{C}$ be a Grothendieck category with a projective generator $G$ and assume that $G$ is small [{\it i.e.,} $Hom(G,{}_{-})$ commutes with all direct sums]. Then $\mathcal{C} \simeq Mod(A^{op})$, for $A = End(G)$. 
\end{thm}

Note that the {\it Gabriel Popescu Theorem} [\ref{Gabriel_Popescu}] gave just a fully faithful embedding with an exact left adjoint and not an equivalence.

\vspace{2mm}

      Now we are ready to discuss a model of \nc projective geometry after
      {\it Artin} and {\it Zhang} \cite{AZ}. However, we would also like to bring into the notice of the readers the works done by {\it Verevkin} (see \cite{Ver}). We have invested enough time in
      trying to convince ourselves that the categorical language should be
      adopted. We shall make no exception to that leitmotif here, but the
      rules of the game will be changed slightly. Fix an algebraically closed
      field $k$; then all categories will be assumed to be {\it $k$-linear} abelian categories [{\it i.e.,}
      the bifunctor $Hom$ ends up in $Mod(k)$]. Since in commutative algebraic geometry one mostly deals with finitely generated $k$-algebras, which are noetherian, here we assume that our $k$-algebras are at least {\it right noetherian}. Let $R$ be a graded algebra. Then we introduce some more categories (which are all $k$-linear):

\begin{eqnarray*}
\begin{split}
Gr(R)   :=\; &\text{category of $\ZZ$-graded right $R$-modules, with degree $0$ morphisms.} \\
Tor(R)  :=\; &\text{full subcategory of $Gr(R)$ generated by torsion modules}\\
           & \text{({\it i.e.,} $M$ such that $\forall$ $x\in M$,
  $xR_{\geqslant s} = 0$ for some $s$), which is thick.} \\
QGr(R)  :=\; &\text{the quotient category $Gr(R)/Tor(R)$ (refer to the sub-section on} \\
             &\text{quotient construction in the previous chapter).} 
\end{split}
\end{eqnarray*}

\begin{rem} \label{notation}
{\bf Standard Convention}. If $XYuvw(\dots)$ denotes an abelian category, then
we shall denote by $xyuvw(\dots)$ the full sub-category consisting of
noetherian objects and if $A,B,\dots,M,N,\dots$ denote objects in $Gr(R)$ then
we shall denote by $\mathcal{A},\mathcal{B},\dots,\mathcal{M},\mathcal{N},\dots$ the corresponding objects in $QGr(R)$.

Some people denote $QGr(R)$ by $Tails(R)$, but we shall stick to our notation. We denote the quotient functor $Gr(R)\map QGr(R)$ by $\pi$. It has a right adjoint functor $\omega : QGr(R)\map Gr(R)$ and so, for all $M\in Gr(R)$ and $\Sheaf{F}\in QGr(R)$ one obtains

\beqn
Hom_{QGr(R)}(\pi M,\Sheaf{F}) \cong Hom_{Gr(R)}(N,\omega\Sheaf{F}).
\eeqn

The $Hom$'s of $QGr(R)$ take a more intelligible form with the assumptions on
$R$. It turns out that for any $N\in gr(R)$ and $M\in Gr(R)$

\beq \label{Hom_desc}
Hom_{QGr(R)}(\pi N,\pi M) \cong \underset{\rightarrow}{Lim}
Hom_{Gr(R)}(N_{\geqslant n},M)
\eeq

For any functor $F$ from a $k$-linear category $\mathcal{C}$ equipped with an autoequivalence $s$, we denote by $\underline{F}$ the graded analogue of $F$ given by $\underline{F}(A):=\underset{n\in\ZZ}{\oplus}F(s^nA)$ for any $A\in Ob(\mathcal{C})$. Further, to simplify notation we shall sometimes denote $s^nA$ by $A[n]$ when there is no chance of a confusion.
\end{rem}

\noindent
Keeping in mind the notations introduced above we have

\begin{lem}
$\omega\pi M \cong \underset{\rightarrow}{Lim} \ul{Hom}_R (R_{\geqslant n}, M)$
\end{lem}

\noindent
{\bf Sketch of proof:}

\beqn
\begin{split}
\omega\pi M & = \ul{Hom}_R (R,\omega\pi M) \hspace{3mm} \text{[since $R\in
  gr(R)$]} \\
            & = \underset{d\in \ZZ}{\oplus} Hom_{Gr(R)} (R,\omega\pi M[d]) \\
            & = \underset{d\in \ZZ}{\oplus} Hom_{QGr(R)}(\pi R,\pi M[d])
  \hspace{3mm} \text{[by adjointness of $\pi$ and $\omega$]} \\
            & = \underset{d\in \ZZ}{\oplus} \underset{\rightarrow}{Lim}
            Hom_{Gr(R)} (R_{\geqslant n},M[d]) \hspace{3mm} \text{[by
            [\ref{Hom_desc}]]} \\
            & = \underset{\rightarrow}{Lim} \underset{d\in \ZZ}{\oplus}
            Hom_{Gr(R)} (R_{\geqslant n}, M[d]) \\
            & = \underset{\rightarrow}{Lim} \ul{Hom}_R (R_{\geqslant n},M)
\end{split}
\eeqn

\QED

The upshot of this lemma is that, there is a natural equivalence of functors $\omega \simeq \ul{Hom} (\mathcal{A},{}_{-})$.

\vspace{2mm}
\begin{center}
\underline{\bf $Proj\, R$}
\end{center}

\vspace{3mm}

Let $X$ be a projective scheme with a line bundle $\Sheaf{L}$. Then the homogeneous coordinate ring $B$ associated to $(X,\Sheaf{L})$ is defined by the formula $B = \underset{n\in \mathbb{N}}{\oplus}\Gamma(X,\Sheaf{L}^{n})$ with the obvious multiplication. Similarly, if $\Sheaf{M}$ is a quasi-coherent sheaf on $X$, $\Gamma_h(\Sheaf{M})=\underset{n\in \mathbb{N}}{\oplus}\Gamma(X,\Sheaf{M}\otimes\Sheaf{L}^{n})$ defines a graded $B$-module. Thus, the compostion of $\Gamma_h$ with the natural projection from $Gr(B)$ to $QGr(B)$ yields a functor ${\bar{\Gamma}}_h:Qcoh(X)\map QGr(B)$. This functor works particularly well when $\Sheaf{L}$ is ample as is evident from the following fundamental result due to {\it Serre}. 

\begin{thm} \cite{Ser}
1. Let $\Sheaf{L}$ be an ample line bundle on a projective scheme $X$. Then the functor ${\bar{\Gamma}}_h({}_{-})$ defines an equivalence of categories between $Qcoh(X)$ and $QGr(B)$. 

2. Conversely, if $R$ is a commutative connected graded $k$-algebra, that is $R_0 = k$ and it is generated by $R_1$ as an $R_0$-algebra, then there exists a line bundle $\Sheaf{L}$ over $X = Proj(R)$ such that $R = B(X,\Sheaf{L})$, up to a finite dimensional vector space. Once again, $QGr(R)\simeq Qcoh(X)$. 
\end{thm}

In commutative algebraic geometry one defines the $Proj$ of a graded ring to be the set of all homogeneous prime ideals which do not contain the augmentation ideal. This notion is not practicable over arbitrary algebras. However, {\it Serre's} theorem filters out the essential ingredients to define the $Proj$ of an arbitrary algebra. The equivalence is controlled by the category $Qcoh(X)$, the structure sheaf $\sheaf_X$ and the autoequivalence given by tensoring with $\Sheaf{L}$, which alludes to the polarization of $X$. Borrowing this idea we get to the definition of $Proj$. Actually one should have worked with a $\ZZ$-graded algebra $R$ and defined its $Proj$ but it has been shown in \cite{AZ} that, with the definition to be provided below, $Proj\, R$ is the same as $Proj\, R_{\geqslant 0}$. Hence, we assume that $R$ is an $\mathbb{N}$-graded $k$-algebra. $Gr(R)$ has a shift operator $s$ such that $s(M)=M[1]$ and a special object, $R_R$. We can actually recover $R$ from the assortment $(Gr(R),R_R,s)$ by

\begin{eqnarray*}
R = \underset{i\in\mathbb{N}}{\oplus} Hom(R_R,s^{i}(R_R))
\end{eqnarray*}

and the composition is given as follows: $a\in R_i$ and $b\in R_j$, then $ab = s^{j}(a)\circ b\in R_{i+j}$. 

 Let $\mathcal{R}$ denote the image of $R$ in $QGr(R)$ and we continue to denote by $s$ the autoequivalence induced by $s$ on $QGr(R)$.

\begin{defn} ($Proj\, R$)

The triple $(QGr(R),\mathcal{R},s)$ is called the projective scheme of $R$ and is denoted $Proj\,R$. Keeping in mind our convention we denote $(qgr(R),\mathcal{R},s)$ by $proj\,R$. This is also equally good due to [\ref{noeVsLocNoe}]. 
\end{defn}

\begin{center}
\underline{\bf Characterization of $Proj\, R$}
\end{center}

\vspace{3mm}

Having transformed {\it Serre's} theorem into a definiton, it is time to address the most natural question - which triples $(\mathcal{C},\mathcal{A},s)$ are of the form $Proj(R)$ for some graded algebra $R$? This problem of characterization has been dealt with comprehensively by {\it Artin} and {\it Zhang}. We would be content by just taking a quick look at the important points. Let us just acquaint ourselves with morphisms of such triples. A morphism between $(\mathcal{C},\mathcal{A},s)$ and $(\mathcal{C}',\mathcal{A}',s')$ is given by a triple $(f,\theta,\mu)$, where $f:\mathcal{C}\map\mathcal{C}'$ is a $k$-linear functor, $\theta:f(\mathcal{A})\map \mathcal{A}'$ is an isomorphism in $\mathcal{C}'$ and $\mu$ is a natural isomorphism of functors $f\circ s\map s'\circ f$. The question of characterization is easier to deal with when $s$ is actually an automorphism of $\mathcal{C}$. To circumvent this problem, an elegant construction has been provided in \cite{AZ} whereby one can pass on to a different triple, where $s$ becomes necessarily an automorphism. If $s$ is an automorphism one can take negative powers of $s$ as well and it becomes easier to define the graded analogues of all functors (refer to [\ref{notation}]). Sweeping that discussion under the carpet, henceforth, we tacitly assume that $s$ is an automorphism of $\mathcal{C}$ (even though we may write $s$ to be an autoequivalence). 

 The definition of $Proj$ was conjured up from {\it Serre's} theorem where the triple was $(Qcoh(X),\sheaf_X ,{}_{-}\otimes\Sheaf{L})$. Of course, one can easily associate a graded algebra to $(\mathcal{C},\mathcal{A},s)$. 

\begin{eqnarray*}
\Gamma_h(\mathcal{C},\mathcal{A},s) = \underset{n\geqslant 0}{\oplus} Hom(\mathcal{A},s^n\mathcal{A}) 
\end{eqnarray*}

with multiplication $a.b = s^n(a)b$ for $a\in Hom(\mathcal{A},s^m\mathcal{A})$ and $b\in Hom(\mathcal{A},s^n\mathcal{A})$. 

But $\Sheaf{L}$ has to be ample and we need a notion of ampleness in the categorical set-up. 

\begin{defn} (Ampleness)

Assume that $\mathcal{C}$ is locally noetherian. Let $\mathcal{A}\in Ob(\mathcal{C})$ be a noetherian object and let $s$ be an autoequivalence of $\mathcal{C}$. Then the pair $(\mathcal{A},s)$  is called ample if the following conditions hold:

1. For every noetherian object $\mathcal{O}\in Ob(\mathcal{C})$ there are positive integers $l_1,\dots,l_p$ and an epimorphism from $\overset{p}{\underset{i=0}{\oplus}}\mathcal{A}(-l_i)$ to $\mathcal{O}$.

2. For every epimorphism between noetherian objects $\mathcal{P}\map\mathcal{Q}$ the induced map $Hom(\mathcal{A}(-n),\mathcal{P})\map Hom(\mathcal{A}(-n),\mathcal{Q})$ is surjective for $n\gg 0$.
\end{defn}

\begin{rem}
The first part of this definition corresponds to the standard definition of an ample sheaf and the second part, to the homological one. 
\end{rem}

Now we are in good shape to state one part of the theorem of {\it Artin} and {\it Zhang} which generalizes that of {\it Serre}.

\begin{thm} \label{characterization}
Let $(\mathcal{C},\mathcal{A},s)$ be a triple as above such that the following conditions hold:

(H1) $\mathcal{A}$ is noetherian,

(H2) $A := Hom(\mathcal{A},\mathcal{A})$ is a right noetherian ring and $Hom(\mathcal{A},\mathcal{M})$ is a finite $A$-module for all noetherian $\mathcal{M}$, and

(H3) $(\mathcal{A},s)$ is ample.

Then $\mathcal{C}\simeq QGr(B)$ for $B=\Gamma_h (\mathcal{C},\mathcal{A},s)$. Besides, $B$ is right noetherian. 
\end{thm}

The converse to this theorem requires an extra hypothesis, which is the so-called $\chi_1$ condition. One could suspect, and rightly so, that there is a $\chi_n$ condition for every $n$. They are all some kind of condition on the graded $Ext$ functor. However, they all look quite mysterious at a first glance. Actually most naturally occurring algebras satisfy them but the reason behind their occurrence is not well understood. We shall discuss them in some cases later but we state a small proposition first.

\begin{prop} \label{chi_equiv}

Let $M\in Gr(B)$ and fix $i\geqslant 0$. There is a right $B$-module structure on $\underline{Ext}^n_B (B/B_+,M)$ coming from the right $B$-module structure of $B/B_+$. Then the following are equivalent:

1. for all $j\leqslant i$, $\ul{Ext}^j_B (B/B_{+},M)$ is a finite $B$-module;

2. for all $j\leqslant i$, $\ul{Ext}^j_B (B/B_{\geqslant n},M)$ is finite for all $n$; 

3. for all $j\leqslant i$ and all $N\in Gr(B)$, $\ul{Ext}^j_B (N/N_{\geqslant n},M)$ has a right bound independent of $n$; 

4. for all $j\leqslant i$ and all $N\in Gr(B)$, $\underset{\rightarrow}{Lim}\ul{Ext}^j_B (N/N_{\geqslant n},M)$ is right bounded.
\end{prop}

       The proof is a matter of unwinding the definitions of the terms suitably and then playing with them. We shall do something smarter instead - make a definition out of it.

\begin{defn} ($\chi$ conditions)

A graded algebra $B$ satisfies $\chi_n$ if, for any finitely generated graded $B$-module $M$, one of the equivalent conditions of the above proposition is satisfied (after substituting $i = n$ in them). Moreover, we say that $B$ satisfies $\chi$ if it satisfies $\chi_n$ for every $n$. 
\end{defn}

\begin{rem} 

Since $B/B_+$ is a finitely generated $B_0$-module ($B_0 = k$) we could have equally well required the finiteness of $\underline{Ext}^n_B (B/B_+,M)$ over $B_0 = k$ {\it i.e.,} $dim_k \underline{Ext}^n_B (B/B_+,M) < \infty$ for $\chi_n$. 
\end{rem}

Let $B$ be an $\mathbb{N}$-graded right noetherian algebra and $\pi : Gr(B)\map QGr(B)$. 

\begin{thm}
If $B$ satisfies $\chi_1$ as well, then (H1), (H2) and (H3) hold for the triple $(qgr(B),\pi B,s)$. Moreover, if $A=\Gamma_h (QGr(B),\pi B,s)$, then $Proj\,B$ is isomorphic to $Proj\,A$ via a canonical homomorphism $B\map A$. [We have a canonical map $B_n = Hom_B(B,B[n])\map Hom(\pi B,\pi B[n]) = A_n$ given by the functor $\pi$.]
\end{thm}

The proofs of the these theorems are once again quite long and involved. So they are left out. What we need now is a good cohomology theory for studying such \nc projective schemes. 

\begin{center}
\underline{\bf Cohomology of $Proj\, R$}
\end{center}

\vspace{3mm}

The following rather edifying theorem due to {\it Serre} gives us some insight into the cohomology of projective (commutative) spaces. 

\begin{thm} \cite{Har}
Let $X$ be a projective scheme over a noetherian ring $A$, and let $\sheaf_X(1)$ be a very ample invertible sheaf on $X$ over $Spec\,A$. Let $\Sheaf{F}$ be a coherent sheaf on $X$. Then:

1. for each $i\geqslant 0$, $H^i(X,\Sheaf{F})$ is a finitely generated $A$-module,

2. there is an integer $n_0$, depending on $\Sheaf{F}$, such that for each $i>0$ and each $n\geqslant n_0$, $H^i(X,\Sheaf{F}(n))=0$. 
\end{thm}

There is an analogue of the above result and we zero in on that. We have already come across the $\chi$ conditions, which have many desirable consequences. Actually the categorical notion of ampleness doesn't quite suffice. For the desired result to go through, we need the algebra to satisfy $\chi$ too. Without inundating our minds with all the details of $\chi$ we propose to get to the point {\it i.e.,} cohomology. Set $\pi R =\mathcal{R}$. On a projective (commutative) scheme $X$ one can define the sheaf cohomology of $\Sheaf{F}\in Coh(X)$ as the right derived functor of the global sections functor {\it i.e.,} $\Gamma$. But $\Gamma(X,\Sheaf{F})\cong Hom_{\sheaf_X}(\sheaf_X,\Sheaf{F})$. Buoyed by this fact and having faith in the people working in ``motivic'' areas, who dream of constructing a universal cohomology theory for varieties (see \cite{FSV}) via the $Ext$ groups, we try to define the cohomology for every $\mathcal{M}\in qgr(R)$ as

\beqn
H^n(\mathcal{M}) := Ext^n_{\mathcal{R}}(\mathcal{R},\mathcal{M})
\eeqn

 However, taking into consideration the graded nature of our objects we also define the following:

\begin{eqnarray*}
\underline{H}^n(Proj\,R,\mathcal{M}) := \ul{Ext}^n_{\mathcal{R}}(\mathcal{R},\mathcal{M}) = \underset{i\in\ZZ}{\oplus} Ext^n_{\mathcal{R}}(\mathcal{R},\mathcal{M}[i]) 
\end{eqnarray*}

The category $QGr(R)$ has enough injectives and one can choose a nice ``minimal'' injective resolution of $\mathcal{M}$ to compute its cohomologies, the details of which are available in the chapter 7 of \cite{AZ}.   

Let $M\in Gr(R)$ and write $\mathcal{M} = \pi M$. Then one should observe that

\beq    \label{cohom_desc}
\begin{split}
\ul{H}^n(\mathcal{M}) & = \ul{Ext}^n_{QGr(R)}(\mathcal{R},\mathcal{M}) \\
                      & \cong \underset{\rightarrow}{Lim}\ul{Ext}^n_{Gr(R)}
                      (R_{\geqslant n},M) \hspace{3mm} [by [\ref{Hom_desc}]] 
\end{split}
\eeq

\vspace{2mm}
As $R$-modules we have the following exact sequence,

\begin{eqnarray*}
0\map R_{\geqslant n}\map R\map R/R_{\geqslant n}\map 0
\end{eqnarray*}

For any $M\in Gr(R)$, the associated $\underline{Ext}$ sequence in $Gr(R)$ looks like

\begin{eqnarray*}
\dots\underline{Ext}^j(R/R_{\geqslant n},M)\map\underline{Ext}^j(R,M)\map\underline{Ext}^j(R_{\geqslant n},M)\map\dots
\end{eqnarray*}

Since $R$ is projective as an $R$ module, $\underline{Ext}^j(R,M)=0$ for every
$j\geqslant 1$. Thus, we get the following exact sequence

\beq \label{exactseq}
\hspace{10mm} 0\rightarrow\underline{Hom}(R/R_{\geqslant n},M)\rightarrow M\rightarrow\underline{Hom}(R_{\geqslant n},M)\rightarrow\underline{Ext}^1(R/R_{\geqslant n},M)\rightarrow 0
\eeq

and, for every $j\geqslant 1$, an isomorphism

\beq \label{Ext_relation}
\underline{Ext}^j(R_{\geqslant n},M)\cong\underline{Ext}^{j+1}(R/R_{\geqslant n},M)
\eeq

The following theorem is an apt culmination of all our efforts. 

\begin{thm}  (Serre's finiteness theorem)

Let $R$ be a right noetherian $\mathbb{N}$-graded algebra satisfying $\chi$, and let $\Sheaf{F}\in qgr(R)$. Then,

(H4) for every $j\geqslant 0$, ${H}^j(\Sheaf{F})$ is a finite right $R_0$-module, and 

(H5) for every $j\geqslant 1$, $\ul{H}^j(\Sheaf{F})$ is right bounded; {\it i.e.,} for $d \gg 0$, ${H}^j(\Sheaf{F}[d])=0$. 
\end{thm}

\noindent
{\bf Sketch of proof:}

Write $\Sheaf{F} = \pi M$ for some $M\in gr(R)$. Suppose that $j=0$. Since $\chi_1(M)$ holds, $\ul{Ext}^i_R (R/R_{\geqslant n},M)$ is a finite $R$-module for each $i=1,2$ and together with [\ref{exactseq}] it implies that $\omega\Sheaf{F} \cong \ul{H}^0(\Sheaf{F})$ is finite (recall $\omega$ from [\ref{notation}]). Now taking the $0$-graded part on both sides we get {\it (1)} for $j=0$.    

Suppose that $j\geqslant 1$. Since $R$ satisfies $\chi_{j+1}$, invoking proposition [\ref{chi_equiv}] we get 

\beqn
\underset{\rightarrow}{Lim} \ul{Ext}^{j+1}_R (R/R_{\geqslant n},M)
\eeqn

is right bounded. Combining [\ref{cohom_desc}] and [\ref{Ext_relation}] this
equals $\ul{H}^j (\Sheaf{F})$. This immediately proves {\it (2)} as
$\ul{H}^j(\Sheaf{F})_d = H^j(\Sheaf{F}[d])$. We now need left boundedness and local finiteness of $\ul{H}^j(\Sheaf{F})$ to finish the proof of {\it (1)} for $j\geqslant 1$. These we have already observed (at least tacitly) but one can verify them by writing down a resolution of $R/R_{\geqslant n}$ involving finite sums of shifts of $R$, and then realizing the cohomologies as sub-quotients of a complex of modules of the form $\ul{Hom}_R (\underset{i=0}{\overset{p}{\oplus}} R[l_i],M)$.    

\QED

Our discussion does not quite look complete unless we investigate the question of the ``dimension'' of the objects that we have defined.

\begin{center}
\underline{\bf Dimension of $Proj\, R$}
\end{center}

\vspace{3mm}

The {\it cohomological dimension} of $Proj\, R$ denoted by $cd(Proj\, R)$ is defined to be 

\begin{eqnarray*}
cd(Proj\, R) := 
  \begin{cases}
     &\text{$sup\{i\; |\; {H}^i(\mathcal{M})\neq 0$ for some $\mathcal{M}\in qgr(R)$\} if it is finite,} \\
     &\text{$\infty$ otherwise.}
  \end{cases}  
\end{eqnarray*} 

\begin{rem}
As ${H}^i$ commutes with direct limits one could have used $QGr(R)$ in the definition of the cohomological dimension.
\end{rem}

The following proposition gives us what we expect from a $Proj$ construction regarding dimension and also provides a useful way of calculating it. 

\begin{prop}

1. If $cd(Proj\, R)$ is finite, then it is equal to $sup\{i\; |\; \underline{H}^i(R)\neq 0\}$. 
\noindent

2. If the left global dimension of $R$ is $d<\infty$, then $cd(Proj\, R)\leqslant d-1$.
\end{prop}

\noindent
{\bf Sketch of proof:}

1. Let $d$ be the cohomological dimension of $Proj\, R$. It is obvious that $sup\{i\; |\; \underline{H}^i(R)\neq 0\}\leqslant d$. We need to prove the other inequality. So we choose an object for which the supremum is attained {\it i.e.,} $\mathcal{M}\in qgr(R)$ such that $H^d(\mathcal{M})\neq 0$ and, hence, $\underline{H}^d(\mathcal{M})\neq 0$. By the ampleness condition we may write down the following exact sequence

\begin{eqnarray*}
0\map \mathcal{N}\map\underset{i=0}{\overset{p}{\oplus}} R[-l_i]\map \mathcal{M}\map 0
\end{eqnarray*}

for some $\mathcal{N}\in qgr(R)$. By the long exact sequence of derived functors $\underline{H}^i$ we have

\begin{eqnarray*}
\dots\map\underset{i=0}{\overset{p}{\oplus}}\underline{H}^d(R[-l_i])\map\underline{H}^d(\mathcal{M})\map\underline{H}^{d+1}(\mathcal{N})=0
\end{eqnarray*}

This says that $\underline{H}^d(R[-l_i])\neq 0$ for some $i$ and hence, $\underline{H}^d(R)\neq 0$. 

\vspace{1mm}
2. It has already been observed that $\underline{H}^i(M)\cong \underset{n\rightarrow \infty}{Lim}\underline{Ext}^i(R_{\geqslant n},M)$ for all $i\geqslant 0$. Now, if the left global dimension of $R$ is $d$, then $\underline{Ext}^j(N,M)=0$ for all $j>d$ and all $N,M\in Gr(R)$. Putting $N=R/R_{\geqslant}$ and using [\ref{Ext_relation}] we get $\underline{H}^d(M)=0$ for all $M\in Gr(R)$. Therefore, $cd(Proj\, R)\leqslant d-1$.   

\QED

\begin{rem}
If $R$ is a noetherian {\it AS-regular} graded algebra, then the {\it Gorenstein} condition can be used to prove that $cd(Proj\, R)$ is actually equal to $d-1$. 
\end{rem}

\begin{center}
\underline{\bf Some Examples} (mostly borrowed from \cite{AZ})
\end{center}

\vspace{3mm}

\begin{ex} {\it (Twisted graded rings)}

Let $\sigma$ be an automorphism of a graded algebra $A$. Then we define a new multiplication $\ast$ on the underlying graded $k$-module $A = \underset{n}{\oplus} A_n$ by 

\begin{eqnarray*}
a\ast b = a\sigma^n(b)
\end{eqnarray*}

where $a$ and $b$ are homogeneous elements in $A$ and $deg(a)=n$. Then algebra is called the {\it twist} of $A$ by $\sigma$ and it is denoted by $A^\sigma$. By \cite{ATV2} and \cite{Zha} $gr(A)\simeq gr(A^\sigma)$ and hence, $proj(A)\simeq proj(A^\sigma)$. 

For example, if $A=k[x,y]$ where $deg(x)=deg(y)=1$, then any linear operator on the space $A_1$ defines an automorphism, and hence, a twist of $A$. If $k$ is an algebraically closed field then, after a suitable linear change of variables, a twist can be brought into one of the forms $k_q[x,y]:=k\{x,y\}/(yx-qxy)$ form some $q\in k$, or $k_j[x,y]:=k\{x,y\}/(x^2+xy-yx)$. Hence, $proj\, k[x,y]\simeq proj\, k_q[x,y]\simeq proj\, k_j[x,y]$. The projective scheme associated to any one of these algebras is the projective line $\proj^1$.
\end{ex}

\begin{ex} {\it (Changing the structure sheaf)}

Though the structure sheaf is a part of the definition of $Proj$, one might ask, given a $k$-linear abelian category $\mathcal{C}$, which objects $\mathcal{A}$ could serve the purpose of the structure sheaf. In other words, for which $\mathcal{A}$ do the conditions (H1), (H2) and (H3) of Theorem [\ref{characterization}] hold? Since (H3) involves both the structure sheaf and the polarization $s$, the answer may depend on $s$. We propose to illustrate the possibilities by the simple example in which $\mathcal{C}=Mod(R)$ when $R=k_1\oplus k_2$, where $k_i = k$ for $i=1,2$ and where $s$ is the automorphism which interchanges the two factors. The objects of $\mathcal{C}$ have the form $V\simeq k_1^{n_1}\oplus k_2^{n_2}$, and the only requirement for (H1), (H2) and (H3) is that both $r_1$ and $r_2$ be not zero simultaneously. 

We have $s^n(V)=k_1^{r_2}\oplus k_2^{r_1}$ if $n$ is odd and $s^n(V)=V$ otherwise. Thus, if we set $\mathcal{A}=V$ and $A=\Gamma_h(\mathcal{C},\mathcal{A},s)$, then $A_n\simeq k_1^{r_1\times r_1}\oplus k_2^{r_2\times r_2}$ if $n$ is even, and $A_n\simeq k_1^{r_1\times r_2}\oplus k_2^{r_2\times r_1}$ otherwise. For example, if $r_1 =1$ and $r_2 =0$, then $A\simeq k[y]$, where $y$ is an element of degree $2$. Both of the integers $r_i$ would need to be positive if $s$ were the identity functor.  
\end{ex}

\begin{ex} {\it (Commutative noetherian algebras satisfy $\chi$)}

Let $A$ be a commutative noetherian $k$-algebra. Then the module structure on $\underline{Ext}^n_A (A/A_+,M)$ can be obtained both from the right $A$-module structure of $A/A_+$ and that of $M$. Choose a free resolution of $A/A_+$, consisting of finitely generated free modules. The cohomology of this complex of finitely generated $A$-modules is given by the $\underline{Ext}$'s, whence they are finite.  

\end{ex}

\begin{ex} {\it (Noetherian AS-regular algebras satisfy $\chi$)}

If $A$ is a noetherian connected $\mathbb{N}$-graded algebra having global dimension $1$, then $A$ is isomorphic to $k[x]$, where $deg(x)=n$ for some $n>0$, which satisfies $\chi$ by virtue of the previous example. In higher dimensions we have the following proposition.

\begin{prop}
Let $A$ be a noetherian AS-regular graded algebra of dimension $d\geqslant 2$ over a field $k$. Then $A$ satisfies the condition $\chi$. 
\end{prop}

\noindent
{\bf Sketch of proof:}

$A$ is noetherian and locally finite (due to finite $GKdim$). For such an $A$ it is easy to check that $\underline{Ext}^j(N,M)$ is a locally finite $k$-module whenever $N,M$ are finite. $A_0$ is finite and hence, $\underline{Ext}^j(A_0 ,M)$ is locally finite for every finite $M$ and every $j$. Since $A$ is connected graded $A_0 = k$. For any $n$ and any finite $A$-module $M$ we first show that $\underline{Ext}^n(A_0,M) = \underline{Ext}^n(k,M)$ is bounded using induction on the projective dimension of $M$. If $pd(M)=0$, then $M=\underset{i=0}{\overset{p}{\oplus}} A[-l_i]$. By the {\it Gorenstein} condition (see definition [\ref{Gorenstein}]) of AS-regular algebra $A$, $\underline{Ext}^n(k,A[-l_i])$ is bounded for each $i$. Therefore, so is $\underline{Ext}^n(k,M)$. If $pd(M)>0$, we choose an exact sequence 

\begin{eqnarray*}
0\map N\map P\map M\map 0
\end{eqnarray*}

where $P$ is projective. Then $pd(N)=pd(M)-1$. By induction, $\underline{Ext}^n(k,N)$ and $\underline{Ext}^n(k,P)$ are bounded, hence, so is $\underline{Ext}^n(k,M)$. Now $A/A_+$ is finite and we have just shown $\underline{Ext}^n(k,M)$ is finite (since bounded together with locally finite implies finite); then $Hom(A/A_+,\underline{Ext}^n(k,M))\cong \underline{Ext}^n(A/A_+,M)$ is locally finite and clearly bounded. Therefore, $\underline{Ext}^n(A/A_+,M)$ is finite for every $n$ and every finite $M$.

\QED 

\begin{rem}
In particular the {\it AS-regular} algebras of dimension $3$ (they are automatically noetherian),  which we have seen in chapter 1, satisfy $\chi$. Also it should be mentioned that all noetherian $\mathbb{N}$-graded $PI$ algebras satisfy $\chi$.
\end{rem}

\end{ex}

There is a host of other examples on algebras satisfying $\chi$ up to varying degrees for which we refer the interested readers to \cite{AZ}, \cite{Rog} and \cite{StZ1}.

\begin{center}
{\bf \underline{Non-commutative smooth proper varieties}}
\end{center}

\vspace{3mm}

The following excerpt is just a leaf out of the article \cite{StV}. Instead of trying to define everything starting from a graded algebra, we now endeavour to define things abstractly and then obtain a description of them in terms of some suitable graded algebras. We first focus on the homological properties that smooth proper varieties have. 

Let $X$ be a smooth, proper and connected $k$-scheme of dimension $n$ and let $\mathcal{C}=Coh(X)$. Then the following properties stare at us discernibly.

\vspace{2mm}
\noindent
{\bf (C1)} $\mathcal{C}$ is noetherian.

\noindent
{\bf (C2)} $\mathcal{C}$ is {\it Ext-finite} [{\it i.e.,} $Ext^i(A,B)<\infty$ for all $A,B\in Ob(\mathcal{C})$ and all $i$].

\noindent
{\bf (C3)} $\mathcal{C}$ has {\it homological dimension} $n$. 

\vspace{2mm}
Unfortunately, easily obtained conditions are not always sufficient and we need to look for more. A fundamental, but slightly more subtle, property of smooth proper schemes is the {\it Serre Duality Theorem} - there exists a dualizing sheaf $\omega$ such that for every $\Sheaf{F}\in\mathcal{C}$ there are natural isomorphisms

\begin{eqnarray*}
H^i(X,\Sheaf{F})\cong Ext^{n-i}(\Sheaf{F},\omega)^\ast \hspace{5mm} \text{[where $\ast$ refers to $k$-dual].}
\end{eqnarray*}

{\it Bondal} and {\it Kapranov} have added to our arsenal a very elegant reformulation of Serre duality \cite{BK}. Let $D^b(\mathcal{G})$ denote the derived category of bounded complexes over an abelian $k$-linear category $\mathcal{G}$. 

\begin{defn} (Serre functor)

A Serre functor on $D^b(\mathcal{G})$ is an autoequivalence $F:D^b(\mathcal{G})\map D^b(\mathcal{G})$ such that there are bifunctorial isomorphisms

\begin{eqnarray} \label{Serre_functor}
Hom(A,B)\cong Hom(B,FA)^\ast
\end{eqnarray} 

which are natural for all $A,B\in Ob(D^b(\mathcal{G}))$.
\end{defn}

Then it is shown that Serre duality can be reinterpreted as saying:

\vspace{1mm}
\noindent
{\bf (C4)} $\mathcal{C}$ satisfies Serre duality in the sense that there exists a Serre functor on $D^b(\mathcal{C})$. 

\vspace{2mm}
We adapt our condition {\bf (C2)} to the graded scenario and demand $\underline{Ext}$-finiteness, instead of $Ext$-finiteness. Classification of abelian $k$-linear categories satisfying {\bf (C1)},{\bf (C2)},{\bf (C3)} and {\bf (C4)} is a tough nut to crack at the moment. However, for hereditary categories, {\it i.e.,} the ones of homological dimension one, the situation is slightly more tractable. 

The most obvious hereditary noetherian $\underline{Ext}$-finite abelian $k$-linear categories satisfying Serre duality in the form of {\bf (C4)} are:

\vspace{2mm}
\noindent
{\bf (E1)} Let X be a smooth projective connected curve over $k$ with function field $K$ and let $\Sheaf{O}$ be a sheaf of hereditary $\sheaf_X$ orders in $M_n(K)$ (see \cite{Rei}). Then one proves exactly as in the commutative case that $\mathcal{C}_1=Coh(\Sheaf{O})$ satisfies Serre duality. The Serre functor is given by tensoring with $Hom(\Sheaf{O},\omega_X)[1]$. 

\vspace{2mm}
\noindent

We now come to three, rather special, examples (all from \cite{StV}) for which {\bf (C1)}-{\bf (C4)} need some work to be verified. 

\vspace{2mm}
\noindent
{\bf (E2)} $\mathcal{C}_2$ consists of the finite dimensional nilpotent representations of the quiver $\tilde{A}_n$ \cite{Rin} (possibly with $n=\infty$) with all arrows oriented in the same direction. The Serre functor is given by rotating one place in the direction of the arrows and shifting one place to the left in the derived category. 

\vspace{2mm}
\noindent
{\bf (E3)} $\mathcal{C}_3=qgr_{\ZZ/2\ZZ}(k[x,y])$ consists of the finitely generated $\ZZ/2\ZZ$-graded $k[x,y]$-modules modulo the finite dimensional ones. The Serre functor is given by the shift $M\longmapsto M(-1,-1)[1]$

\vspace{2mm}
\noindent
{\bf (E4)} The category $\mathcal{C}_3$ has a natural automorphism $\sigma$ of order $2$ which sends a graded module $(M_{ij})_{ij}$ to $(M_{ji})_{ij}$. and which exchanges the $x$ and $y$ action. Let $\mathcal{C}_4$ be the category of $\ZZ/2\ZZ$ equivariant objects of $\mathcal{C}_3$, {\it i.e.,} pairs $(M,\phi)$ where $\phi$ is an isomorphism $M\map \sigma(M)$ satisfying $\sigma(\phi)\phi=id_M$. Properties {\bf (C1)}-{\bf (C4)} follow easily from the fact that they hold for $\mathcal{C}_3$. In particular, the Serre functor is obtained from the one on $\mathcal{C}_3$ in the obvious way. As described, this construction requires $char\; k\neq 2$, but this can be circumvented \cite{RV}.

\vspace{2mm}
\noindent
We now present the most subtle and exotic example.

\vspace{2mm}
\noindent
{\bf (E5)} Let $Q$ be a connected locally finite quiver such that the opposite quiver has no infinite oriented paths. For a vertex $x\in Q$ we have a corresponding projective representation $P_x$ and an injective representation $I_x$ and by our hypotheses $I_x$ is finite dimensional. Let $rep(Q)$ be the finitely presented representations of $Q$. The category of all representations of $Q$ is hereditary and hence, $rep(Q)$ is an abelian category. It is easy to check that $I_x\in rep(Q)$. Hence, the functor $P_x\map I_x$ may be derived to yield an endo-functor $F:D^b(rep(Q))\map D^b(rep(Q))$. This functor $F$ behaves like a Serre functor in the sense that we have natural isomorphisms as in relation [\ref{Serre_functor}], but unfortunately $F$ need not be an autoequivalence. However, there is a formal procedure to quasi-invert $F$ so as to obtain a true Serre functor (see \cite{RV}). This yields a hereditary category $\widetilde{rep}(Q)$ which satisfies Serre duality. Under a technical additional hypothesis (see \cite{RV}) $\mathcal{C}_5 =\widetilde{rep}(Q)$ turns out to be noetherian. 

\vspace{2mm}
We would just like to recall a simple definition here. 

\begin{defn} (connectedness of an abelian category)

An abelian category $\mathcal{G}$ is {\it connected} if, whenever $\mathcal{G}=\mathcal{G}_1\oplus\mathcal{G}_2$, then either $\mathcal{G}_1 =\mathcal{G}$ or $\mathcal{G}_2 = \mathcal{G}$.   
\end{defn}
 
The examples {\bf (E1)}-{\bf (E5)} are all connected and thanks to {\it Reiten} and {\it M. van den Bergh} we know that they are the only hereditary ones satisfying {\bf (C1)},{\bf (C2)} and {\bf (C4)}.

\begin{thm} \cite{RV} \label{classification}
Let $\mathcal{C}$ be a connected noetherian $\underline{Ext}$-finite hereditary category satisfying Serre duality over $k$. Then $\mathcal{C}$ is equivalent to one of the categories {\bf (E1)}-{\bf (E5)}.
\end{thm}

\vspace{2mm}
It would be rather hasty to conclude that we have hit upon the right set of hypotheses for a \nc analogue of a smooth proper variety. In fact, one indication that {\bf (C1)}-{\bf (C4)} may be inadequate is that, for example, they do not distinguish between algebraic and analytic (smooth compact) surfaces. To get around this problem we resort to another property of smooth proper schemes that was also discovered by {\it Bondal} and {\it Kapranov} and is not satisfied by analytic surfaces. 

\begin{defn} (finite type cohomological functor)

A cohomological functor $H:D^b(\mathcal{C})\map mod(k)$ is of finite type if, for any $A\in D^b(\mathcal{C})$,
\beqn
\underset{n}{\sum}dim_k H(A[n]) < \infty .
\eeqn 
\end{defn}

\noindent
Then the appropriate condition is:

\vspace{2mm}
\noindent
{\bf (C5)} Let $\mathcal{C}$ be an $\underline{Ext}$-finite abelian $k$-linear category of finite homological dimension. Then $\mathcal{C}$ is {\it saturated} if every cohomological functor $H:D^b(\mathcal{C})\map mod(k)$ of finite type is representable. [The condition requires that $\mathcal{C}$ be saturated].

\vspace{2mm}
It was shown in \cite{BV} that $Coh(X)$ is saturated when $X$ is a smooth projective scheme and that saturation also hold for categories of the from $mod(\Lambda)$ where $\Lambda$ is a finite dimensional algebra. Further, saturated categories satisfy Serre duality. So it is a stronger criterion. 

Combined with theorem [\ref{classification}], this gives a much more compact classification.

\begin{cor}
Assume that $\mathcal{C}$ is a saturated connected noetherian $\underline{Ext}$-finite hereditary category. Then $\mathcal{C}$ has one of the following forms:

(1) $mod(\Lambda)$ where $\Lambda$ is an indecomposable finite dimensional hereditary algebra (this is a very special case of {\bf (E5)}). 

(2) $Coh(\Sheaf{O})$ where $\Sheaf{O}$ is a sheaf of hereditary $\sheaf_X$ orders (see {\bf (E1)}) over a smooth connected projective curve $X$. 
\end{cor}

\begin{rem}
It can also be shown that the abelian categories occurring in this corollary are of the form $qgr(R)$ for some graded ring $R$ of $GKdim\leqslant 2$. This, in a sense, will prove the classical commutative result that smooth proper curves are projective. 
\end{rem}
  
We are now going to provide some criteria for a module category to be saturated. Interestingly, the seemingly strange $\chi$ conditions show up here as well. For the result stated below, we declare that a connected graded ring $R$ has finite right cohomological dimension, provided the higher right derived functors of the functor $\underset{\map}{Lim}\underline{Hom}_R (R/R_{\geqslant n},{}_{-})$ vanish {\it i.e.,} $R^i\underset{\map}{Lim}\underline{Hom}_R (R/R_{\geqslant n},{}_{-}) = 0$ for $i\gg 0$ [refer to the discussion on the dimension of $Proj\, R$]. 

\begin{thm} \cite{BV}
Let $R$ be a connected graded noetherian ring satisfying the following hypotheses:

1. $R$ and its opposite ring $R^\circ$ satisfy $\chi$;

2. $R$ and $R^\circ$ have finite right cohomological dimension;

3. $qgr(R)$ has finite homological dimension.

\vspace{1mm}
\noindent
Then $qgr(R)$ is saturated. 
\end{thm}

As a justification for bringing in this {\it saturation} criterion we mention the following result.

\begin{prop} \cite{BV}
Assume that $X$ is an analytic $K3$ surface with no curves. Then $Coh(X)$ is not saturated. 
\end{prop}

\vspace{2mm}
With all the rambling above we hope to have augmented our understanding of \nc smooth proper ``curves'' at least. Now it is worth taking a closer look at Serre duality as it only makes sense for categories which have finite homological dimension, whereas in classical algebraic geometry much of the significance of Serre duality lies in its applications to singular varieties, whose categories of coherent sheaves are infinite dimensional. So we investigate a second notion which is applicable more generally.

Suppose that $\mathcal{C}$ is an $\underline{Ext}$-finite noetherian $k$-linear abelian category with a distinguished object $\Sheaf{O}\in Ob(\mathcal{C})$. Then $(\mathcal{C},\Sheaf{O})$ is said to satisfy {\it classical Serre duality} if there exists an object $\omega^\circ\in D^b(\mathcal{C})$ together with a natural isomorphism 

\begin{eqnarray*}
RHom({}_{-},\omega^\circ)\cong RHom(\Sheaf{O},{}_{-})^\ast
\end{eqnarray*}

We can always pass on to $Ind(\mathcal{C})$ to compute these derived functors. Since $\omega^\circ$ represents a functor it is clear that it must be unique if it exists, and, moreover, its existence forces $Hom(\Sheaf{O},{}_{-})$ to have finite cohomological dimension. The following result gives us a class of rings for which $\omega^\circ$ exists. 

\begin{thm} \cite{YZ1} \label{classical_Serre}
Let $(\mathcal{C},\Sheaf{O})=(qgr(R),\pi R)$ for a noetherian connected graded ring $R$. Assume, further, the following:

1. $R$ and its opposite ring $R^\circ$ satisfy $\chi$;

2. $R$ and $R^\circ$ have finite right cohomological dimension. 

\vspace{1mm}
\noindent
Then $(qgr(R),\pi R)$ satisfies classical Serre duality.
\end{thm}

\begin{rem}
If we allow $\omega^\circ$ to lie in $D(QGr(R))$, then its existence does not require the $\chi$ assumption (see \cite{Jor1}). {\it J\o rgensen, Yekutieli} and {\it Zhang} among others have explored the veracity of certain homological results from the commutative purview in this new model under the assumptions of theorem [\ref{classical_Serre}] (see \cite{Jor2}, \cite{Jor3}, \cite{JZ}, \cite{StZ2}, \cite{WZ} and \cite{YZ2}). As a final remark it should be mentioned that {\it Gorenstein} rings of finite left and right injective dimensions satisfy the assumptions of theorem [\ref{classical_Serre}] and hence, classical Serre duality (see \cite{YZ2}). 
\end{rem}

Generally one begins with definitions but we would like to bring the curtain
down on this section with a proposed definition of \nc $\proj^2$ by {\it
  Artin} and {\it Schelter}, which uses bits and pieces of almost everything that we have discussed so far.

\begin{defn} (\nc $\proj^2$) (\cite{StV} definition 11.2.1.) 

A non-commutative $\proj^2$ is a $k$-linear Grothendieck category of the form
$QGr(A)$, where $A$ is an $AS$-regular algebra of dimension $3$, with Hilbert
series $(1-t)^{-3}$. 
\end{defn}

Following the classification of such algebras described in chapter 1, one
gets two different $\proj^2$'s, one of them being the known commutative one and the other, a truly \nc $\proj^2$. 


\newpage

\section{Non-commutative Geometry via Deformation Theory}

\vspace{5mm}

This line of thought was propounded by {\it Laudal} and it provides a different flavour from what we have seen so far. Still being very algebraic in its approach, it manages to develop an interesting ``infinitesimal'' theory, which is also in keeping with the title of the chapter. It is deemed fit that we take a quick look at the commutative deformation theory after {\it Schlessinger} \cite{Sch}.

\vspace{2mm}
Let $k$ be an algebraically closed field. Let $Art_k$ denote the category of {\it Artinian, local} ``commutative'' $k$-algebras having residue field $k$ {\it i.e.,} diagrams of the form 

\xymatrix{
k \ar[r] \ar@{=}[dr] 
      & (A,\mathtt{m}) \ar[d] \\
      & A/\mathtt{m}=k }

We shall denote $Art_k$ by $\mathcal{C}$ and by $\hat{\mathcal{C}}$ we shall understand the pro-completion of $\mathcal{C}$ [it is dual to the construction of $Ind(\mathcal{C})$ (see [\ref{IndC}])]. $\hat{\mathcal{C}}$ has a more intelligible description: its objects are complete local $k$-algebras $R$, such that $R/\mathtt{m}^n \in \mathcal{C}$ for all $n\geqslant 1$, where $\mathtt{m}$ is the unique maximal ideal of $R$.  A functor $F:\mathcal{C}\map Set$ naturally extends to a functor $\hat{F} :\hat{\mathcal{C}}\map Set$ and by some {\it Yoneda} kind of argument we have

\begin{eqnarray*}
\hat{F}(R)\simeq Hom(h_R,F)
\end{eqnarray*}

So $\zeta\in\hat{F}(R)$ induces a morphism $h_R\map F$. 

\begin{defn} (Smoothness of morphisms of functors) \label{smooth}

A morphism $F\map G$ of functors is smooth if for any surjection $B\map A$ in $\mathcal{C}$ the morphism 

\begin{eqnarray*}
F(B)\map F(A)\times_{G(A)} F(B)
\end{eqnarray*} 

is surjective as well. 
\end{defn}

\vspace{2mm}
\noindent
The tangent space of $F$, denoted by $t_F$ $:=$ $F(k[\epsilon])$. 

\begin{defn} (Pro-representing hull)

Let $F:\mathcal{C}\map Set$ be a functor. A couple $(R,\zeta)$, where $R\in\hat{\mathcal{C}}$ and $\zeta\in \hat{F}(R)$, is a pro-representing hull of $F$ if the induced map $h_R\map F$ is smooth and if in addition $t_{h_R}\map t_F$ is a bijection.
\end{defn}

\noindent
Analogously one would say that $(R,\zeta)$ {\it pro-represents} $F$  if the morphism $h_R\map F$ induced by $\zeta$ is an isomorphism.

\noindent
With all these definitions in our fingertips, we can now state the following theorem due to {\it Schlessinger} \cite{Sch}. 

\begin{thm}
Let $F:\mathcal{C}\map Set$ be a covariant functor such that $F(k)$ is singleton. Let $A'\map A$ and $A''\map A$ be morphisms in $\mathcal{C}$ and consider the map,

\begin{eqnarray*}
F(A'\times_A A'') \stackrel{\psi}{\map} F(A')\times_{F(A)} F(A'')
\end{eqnarray*}

Then $F$ has a pro-representing hull if and only if $F$ has the properties {\bf (H1)}, {\bf (H2)} and {\bf (H3)} below:

\noindent
{\bf (H1)} $\psi$ is a surjection whenever $A''\map A$ is a small extension [{\it i.e.,} $ker(A''\map A)$ is a principal ideal $\subset Ann(\mathtt{m}_{A''})$]

\noindent
{\bf (H2)} $\psi$ is a bijection when $A=k$, $A''=k[\epsilon]$. 

\noindent
{\bf (H3)} $dim_k\, t_F < \infty$.

\vspace{1mm}
Further, $F$ is pro-representable by a couple if and only if $F$ satisfies

\noindent
{\bf (H4)} $\psi$ is an isomorphism whenever $A'=A''$ and $A'\map A$ is a small extension. 
\end{thm}  

Let $X\in Sch_k$ and $A\in\mathcal{C}=Art_k$. Then by deformation of $X$ to $A$ we mean a diagram of the form 

\[
\begin{CD}
X @>i>closed\, immersion> Y \\
@VVV @VVflatV \\
Spec\, k @>>> Spec\, A
\end{CD}
\]

such that $X\overset{\sim}{\map} Y\times_{Spec\, A} Spec\, k$. To any morphism
$A\map B$ and $Y$, a deformation of $X$ to $A$, we can associate
$Y\times_{Spec\, A} Spec\, B$, a deformation to $B$. We can package all these
information into one single definition. 

\begin{defn} (Deformation functor $Def_{X|_k}$)

\begin{eqnarray*}
\begin{split}
Def_{X|_k}:\, & Art_k \map Set \\
            & A \longmapsto \text{set of isomorphism classes of deformations
            of $X$ to $A$}
\end{split}
\end{eqnarray*}
\end{defn}

{\it Schlessinger} in \cite{Sch} showed that $Def_{X|_k}$ satisfies {\bf (H1)}
and {\bf (H2)}. For {\bf (H3)} he needed some assumptions on $X$ ({\it e.g.,}
$X$ is proper over $k$). For {\bf (H4)} also he has got a requirement which we are going skip over here. 

\vspace{1mm}
One can easily figure out the right definition of the {\it deformation functor} of other algebraic geometric objects {\it e.g.,} sheaves of modules. With that we come to the end of our cruise through commutative ``deformation theory''. 

\vspace{2mm}
We shall look at deformations which will eventually enable us to enter the
world of \nc algebraic geometry. With some foresight, we shall study the
simultaneous formal deformations of a finite family [the infinite theory can
  be developed suitably as a limiting case of the finite ones] of right
$A$-modules, where $A$ is an associative $k$-algebra. The right place {\it i.e., category} to develop \nc deformation theory turns out to be $\ul{a}_r$, which is a suitable sub-category of $\ul{A}_r$, the category of $r$-pointed $k$-algebras. The objects of $\ul{A}_r$ are the diagrams of $k$-algebras of the following form:

\xymatrix{
k^r \ar[r]^i \ar@{=}[dr] 
      & R \ar[d]^p \\
      & k^r }

A few words about terminology: by a {\it diagram} we shall understand an
oriented graph (could be infinite) with possibly multiple edges between two
vertices. This is essentially a {\it local terminology}. There is no more assumption on a diagram and so such an object may not be a
category ({\it e.g.,} if it does not have any self loop at a vertex, which is
by the way not forbidden, then the identity morphism does not belong to the
diagram for that particular vertex). Since only the map $p$ is needed in most
cases, henceforth we shall also denote the objects of $\ul{A}_r$ by $p:R\map k^r$ or simply by $R$. 

The radical of $R$ is $Rad(R):=ker(p)$. Put
$T_R=(Rad(R)/Rad(R)^2)^\ast$ and call it the {\it tangent space} of $R$. The
category $\ul{a}_r$ is the full sub-category of $\ul{A}_r$ consisting of
objects in which $R$ is Artinian and complete in the $Rad(R)$-adic topology. 

Fix a finite diagram $\mathcal{V}=\{V_1,\dots,V_r\}$ of right $A$-modules. Such a family is called a {\it swarm} if, for all $i,j$

\begin{eqnarray*}
dim_k Ext^1_A(V_i,V_j) < \infty
\end{eqnarray*}

Henceforth we are going to use the terms, ``diagrams'' and ``swarms'', almost interchangeably. For any $R\in \ul{a}_r$ we define a {\it lifting} of $\mathcal{V}$ to $R$ as an $R-A$ bimodule denoted by $\mathcal{V}_R$, together with isomorphisms $\eta_i : k_i\otimes_R\mathcal{V}_R\map V_i$ of right $A$-modules for $1\leqslant i\leqslant r$, such that the left and right $k$-vector space structures coincide and abstractly $\mathcal{V}_R\simeq R\otimes_k \mathcal{V}$ as left $R$-modules. The requirement that $\mathcal{V}_R\simeq R\otimes_k \mathcal{V}$ as left $R$-modules generalizes the flatness condition of commutative deformation theory.  

Let $\mathcal{V}'_R,\mathcal{V}''_R$ be two liftings of $\mathcal{V}$ to $R$. We say that these two liftings are {\it isomorphic} if there exists an isomorphism $\tau:\mathcal{V}'_R\map\mathcal{V}''_R$ of $R-A$-bimodules, such that $\eta''_i\circ(Id\otimes\tau)=\eta'_i$ for all $1\leqslant i\leqslant r$. We refer to these equivalence classes as {\it deformations} of $\mathcal{V}$ to $R$. Although $\eta_i$'s form an integral part of the definition of liftings and deformations, for convenience we shall often suppress them.

\begin{defn} (Non-commutative deformation functor)

It is a covariant functor defined as follows:

\begin{eqnarray*}
\begin{split}
Def_{\mathcal{V}} : &\ul{a}_r\map Set \\
                    & R \longmapsto \text{set of deformations of $\mathcal{V}$ to $R$}
\end{split}
\end{eqnarray*}

Given any morphism $R\map S$ in $\ul{a}_r$ one can easily check that $\mathcal{V}_S = S\otimes_R\mathcal{V}_R$ is a deformation of $\mathcal{V}$ to $S$, which is independent of the lifting chosen. 

\end{defn}

Now it is a theorem due to {\it Laudal}, which illustrates the pro-representability of the \nc deformation functor.

\begin{thm}
The functor $Def_\mathcal{V}$ has a pro-representable hull {\it i.e.,} an object $H(\mathcal{V})=H\in\hat{\ul{a}}_r$, together with a versal family,

\beqn
\tilde{\mathcal{V}} = H\hat{\otimes}_k \mathcal{V}
\eeqn

such that the corresponding morphism of functors in $\hat{\ul{a}}_r$, 

\beqn 
\rho :Hom_{\hat{\ul{a}}_r}(H,{}_{-})\map {Def}_{\mathcal{V}}
\eeqn

is smooth (see [\ref{smooth}]) and an isomorphism at the tangent level [{\it i.e.,} $\rho$ is an isomorphism for all $R\in\ul{a}_r$ for which $Rad(R)^2 = 0$]. 
\end{thm} 

The fairly involved proof of this theorem can be found in \cite{Lau2}. {\it
  Ile} has a proof too, which is actually closer in spirit to that of {\it Schlessinger} in the
commutative case (see \cite{Ile}). One more apposite remark has to be made before we move on.

\begin{rem}  \label{Burnside}
By virtue of the definition of a deformation, $\tilde{V}$ is itself an $H-A$-bimodule. The right $A$-module structure of $\tilde{V}$ defines a homomorphism of $k$-algebras, 

\beqn
\eta : A\map End_H (\tilde{V})=:\sheaf (\mathcal{V})
\eeqn

and the $k$-algebra $\sheaf (\mathcal{V})$ extends the right module action of $A$ on $\mathcal{V}=\{V_1,\dots,V_r\}$. 

\vspace{2mm}
{\bf A generalized Burnside theorem:} Let $A$ be a finite dimensional
$k$-algebra, $k$ being an algebraically closed field. Suppose that the family
$\mathcal{V}=\{V_1,\dots,V_r\}$ contains all simple $A$-modules, then 

\beqn
A\simeq \sheaf (\mathcal{V})
\eeqn

Further, $A$ is Morita equivalent to $H(\mathcal{V})$, which is evident from
the relation between $\sheaf(\mathcal{V})$ and $H(\mathcal{V})$.
\end{rem}

We have barely been able to scratch the surface of \nc deformation theory, and
it is time to leave it at that and plunge into \nc algebraic
geometry. Consider $X:= Spec\, A$, where $A$ is a commutative finite type
$k$-algebra. Then the closed points of $X$ are given by the maximal ideals and
to any such closed point $x\in X$ we can associate the $A$-module $k(x)$ which
is the residue field at $x$. Suppose for $\mathcal{V}$ we take the singleton
set $\{k(x)\}$ for a fixed $x\in X$. It is almost tautological to say that $X$
is the moduli space of its closed points (not all points) and one has to
ruminate a little to be convinced of the fact that $H(\mathcal{V})$, the hull
of the deformation functor $Def_\mathcal{V}$, will be given by
$\hat{A}_{\mathtt{m}_x}$, which is the $\mathtt{m}_x$-adic completion of the
local $k$-algera, $\sheaf_{X,x}=A_{\mathtt{m}_x}$. A section $f$ of the structure sheaf can be read off from its germs $\hat{f}_x\in\hat{A}_{{\mathtt{m}}_x}$ at the different points $x\in X$. This is what we are going to cash in on and the new notion of {\it structure sheaf} will be obtained from germane deformations. 

For \nc rings we have a good cohomology theory, namely {\it
  Hochschild Cohomology} (and also Connes' Cyclic cohomology). In the
  evolutionary process of mathematical objects, rings seem to be the
  ancestors of pre-additive categories and now we shall see the evolution of
  the cohomology theory which will be applicable to small abelian
  categories. By the way, evolution is a dynamic process.

\vspace{1mm}
For any small abelian categore $\mathcal{C}$ we shall define a new category
$Mor\mathcal{C}$. 

1. The objects are just morphisms in $\mathcal{C}$ 

2. If $\phi ,{\phi}'$ are two arrows in $\mathcal{C}$ then
   $Hom_{Mor\mathcal{C}}(\phi ,{\phi}')$ is the set of commutative diagrams 

\[
\begin{CD}
\ast @<\psi<< \ast \\
@V\phi VV @VV\phi' V \\
\ast @>>\psi'> \ast
\end{CD}
\] 

We shall simply write $(\psi,\psi')$ for such an arrow in $Mor\mathcal{C}$.

Now consider any covariant functor $G:Mor\mathcal{C}\map Ab$. We shall define a
chain complex, $D^\ast (\mathcal{C},G)$ as follows:

\beqn
\begin{split}
D^0(\mathcal{C},G) = & \underset{c\in Ob(\mathcal{C})}{\prod} G(Id_c) \\
D^p(\mathcal{C},G) = &
\underset{c_0 \overset{\psi_1}{\rightarrow} c_1\dots\overset{\psi_p}{\rightarrow}c_p}{\prod}
G(\psi_1\circ\dots\circ\psi_p)
\end{split}
\eeqn

where the indices are strings of morphisms $\psi_i:c_{i-1}\map c_i$ of length $p$. The differential
$d^p : D^p(\mathcal{C},G)\map D^{p+1}(\mathcal{C},G)$ is defined in the
obvious manner {\it i.e.,} if $\zeta\in D^p(\mathcal{C},G)$ then the
$\psi_1\circ\dots\circ\psi_{p+1}$th entry of $d^p(\zeta)$, denoted by $d^p(\zeta)_{(\psi_1,\dots,\psi_{p+1})}$ is:

\beqn
\begin{split}
d^p(\zeta)_{(\psi_1,\dots,\psi_{p+1})}= & \, G(\psi_2\circ\dots\circ\psi_{p+1}\rightarrow\psi_1\circ\dots\circ\psi_{p+1})(\zeta_{(\psi_2,\dots,\psi_{p+1})}) \\
& +\underset{i=1}{\overset{p}{\sum}}(-1)^i\zeta_{(\psi_1,\dots,\psi_i\circ\psi_{i+1},\dots,\psi_{p+1})} \\
& + (-1)^{p+1}G(\psi_1\circ\dots\circ\psi_p\rightarrow\psi_1\circ\dots\circ\psi_{p+1})(\zeta(\psi_1,\dots,\psi_p))
\end{split}
\eeqn

where
$(\psi_2\circ\dots\circ\psi_{p+1}\rightarrow\psi_1\circ\dots\circ\psi_{p+1})$
and
$(\psi_1\circ\dots\circ\psi_{p}\rightarrow\psi_1\circ\dots\circ\psi_{p+1})$
are arrows in $Mor\mathcal{C}$.

Now one has to check that $d^{p+1}\circ d^p = 0$ and then we may summarize our
efforts by saying that 

\beqn
\begin{split}
D^\ast :  Funct(Mor\mathcal{C},Ab) &\map \text{Cat. of complexes of ab. gps.} \\
            G &\longmapsto (D^p(\mathcal{C},G),d^p)_{p\geqslant 0}
\end{split}
\eeqn

is a functor, which has been shown to be exact in \cite{Lau1}. Note that, the functor
$D^\ast$ is exact, but given a functor $G:Mor\mathcal{C}\map Ab$, the complex
$(D^p(\mathcal{C},G),d^p)$ can have non-trivial cohomologies. 

We define the cohomology of the category $\mathcal{C}$ with respect to a
functor $G:Mor\mathcal{C}\map Ab$ as 

\beqn
H^\ast(\mathcal{C},G) := H^\ast(D^p(\mathcal{C},G),d^p)
\eeqn

Given an associative $k$-algebra $A$, let $\mathcal{C}$ be any sub-diagram of
$Mod(A)$, the category of right $A$-modules. Let $\pi : \mathcal{C}\map
Mod(k)$ denote the forgetful ``functor''. Then we can manufacture a functor
$Hom_\pi : Mor\mathcal{C}\map Mod(k)$. 

\beqn
\begin{split}
Hom_\pi : Mor\mathcal{C} &\map Mod(k) \\
          c_1\overset{\psi}{\rightarrow} c_2 &\longmapsto
Hom_k(\pi(c_1),\pi(c_2))
\end{split}
\eeqn 

\noindent
This allows us to define a transient object, which is $O_o(\mathcal{C},Hom_\pi):=
  H^0(\mathcal{C},Hom_\pi)$. I am afraid, but we are going to come across a plethora of $O$'s now. We shall simplify notation by writing just
  $O_o(\mathcal{C})$ instead of $O_o(\mathcal{C},Hom_\pi)$. It is clear from the definition of the cohomology that
  $O_o(\mathcal{C})\subset \underset{c\in\mathcal{C}}{\prod}End_k(\pi(c))$ and
  it inherits its algebra structure from $\underset{c\in\mathcal{C}}{\prod}End_k(\pi(c))$. Moreover, it also has a canonical
  projection onto each $End_k(\pi(c))$ making $c$ into an
  $O_o(\mathcal{C})$-module. There is also a canonical homomorphism of
  $k$-algebras $\eta_0 : A\map O_o(\mathcal{C})$ as every $a\in A$ is seen to give rise
  to a $0$-cocycle. 

We had started with $\mathcal{C}$, which was a sub-diagram of $Mod(A)$, and we
have realized that each $c\in \mathcal{C}$ is also an $O_o(\mathcal{C})$-module. So
we may play the same game considering $\mathcal{C}$ as a sub-diagram of
$Mod(O_o(\mathcal{C}))$ and taking a new forgetful functor (restricting to $\mathcal{C}$) 

\beqn 
\pi_0 : Mod(O_o(\mathcal{C}))\map Mod(k)
\eeqn

The striking thing is that $O_o(\mathcal{C},\pi_0)\cong
O_o(\mathcal{C},\pi)$ and that we get nothing new.  
 
\begin{rem}
For a commutative $k$-algebra $A$ of finite type let $Prim(A)$ denote the
sub-category of $Mod(A)$ comprising the indecomposable modules of the form
$A/\mathtt{q}$ for some primary ideal $\mathtt{q}$. Then, if $A=k[\epsilon]$,
it turns out that $O_o(Prim(A),\pi)\simeq$ $\{ 2\times 2$ upper triangular
  matrices with entries in $k\}$. Here, the algebra $A$ has nilpotent
  elements.
However, if $A$ were reduced we would have had an isomorphism via $\eta_0$. 
\end{rem}

\noindent
There is a {\it Jacobson topology} on the set $Prim(A)$ ($A$ reduced) defined as follows:

Let $a\in A$ and consider the subset $D(a)$ of $Prim(A)$ defined by the
objects $M$ for which $a\notin Ann(M)$. Then clearly $D(a)\cap D(b) = D(ab)$
and hence, they form a basis of the topology generated by them. What is more interesting is that, $D(a)$ is
simply the localization of $Prim(A)$ at the multiplicative subset generated by
the powers of $a$ and for each $a$ there is a canonical isomorphism 

\beqn
O_o(D(a),\pi) \simeq A_{(a)}= \sheaf_{Spec\, A}(D(a))
\eeqn

In fact, as schemes $Spec\, A$ and $Prim(A)$ are isomorphic. So {\it Jacobson
  topology} is a reincarnation of Zariski topology over $Prim(A)$. One
could have thought of defining $Prim(A)$'s as \nc schemes but, actually there
  is still room for improvement. 

Recall the $\sheaf(\mathcal{V})$ construction (see
[\ref{Burnside}]) for a swarm $\mathcal{V}$ of $A$-modules. It came equipped
with a canonical $k$-algebra homomorphism $\eta:A \map \sheaf(\mathcal{V})$. 
Also from the very definition of the terms, we obtain a canonical $k$-algebra
homomorphism $\rho_0:\sheaf(\mathcal{V})\map O_o(\mathcal{V})$. 
We summarize the whole state of affairs in the following commutative diagram
for a finite diagram $\mathcal{V}$ of $A$-modules.

\begin{center}

\xymatrix{
  A \ar[r]^\eta \ar[dr]_{\eta_0}
     & \sheaf(\mathcal{V}) \ar[d]^{\rho_0} \\
     & O_o(\mathcal{V}) }

\end{center}

The arrows extend the module actions on $\mathcal{V}$.

\spc

Now we introduce a notation which should have been done much earlier. For any
diagram of $A$-modules $\mathcal{V}:= \{V_1,\dots ,V_r\}$ let $|\mathcal{V}|$ denote the underlying
set of modules {\it i.e.,} the diagram stripped off its morphisms.
So far we have not at all made use of the
  fact that there could be morphisms between the deformed modules of the swarm
  (or the diagram) and so the $\sheaf(\mathcal{V})$ that we have constructed
  is actually $\sheaf(|\mathcal{V}|)$. There is a way to incorporate the
{\it incidence datum} of the swarm of modules into the definition of the algebra $\sheaf(\mathcal{V})$, which we are unfortunately not able to discuss. However, we
just present the definition with a host of strange symbols with an appeal to
the interested readers to look them up in \cite{Lau3}. The nomenclature has a direct relationship with physics.

\begin{defn} (Algebra of pre-observables) 

The $k$-algebra of pre-observables $O(\mathcal{V},\pi)$ of the finite swarm
$\mathcal{V}$, is the sub-algebra of $(H(\mathcal{V})\otimes End_\pi
(\mathcal{V}))$ [generated by the morphisms in $\pi\mathcal{V}$] commuting, for any $S\in \ul{\hat{a}}_r$, via the morphism

\beq   \label{pre-obs}
\kappa_S : (H(\mathcal{V})\otimes End_\pi (\mathcal{V})) \map (S\otimes End_k
(\mathcal{V}))
\eeq

induced by any surjective $k$-algebra homomorphism 

\beqn 
H(\mathcal{V})\map S
\eeqn

with the corresponding representation of the path
algebra of $\mathcal{V}$ (thought of as a quiver) in $(S\otimes End_k(\mathcal{V}))$.
\end{defn}

Now we need to take a look at the problem of functoriality. This $O$-construction is
not functorial with respect to the inclusions of diagrams of modules and their corresponding path algebras. An arrow of the form $k[\mathcal{C}_0]\hookrightarrow k[\mathcal{C}]$
induces an arrow like the one below:

\beqn
o(\mathcal{C}_0\subset\mathcal{C}): (H(\mathcal{C})\otimes End (\mathcal{C}))\map (H(\mathcal{C}_0)\otimes End (\mathcal{C}_0))
\eeqn

We now desire that there be a natural morphism 

\beqn
O(\mathcal{C}\subset\mathcal{C}_0): O(\mathcal{C},\pi)\map
O(\mathcal{C}_0,\pi)
\eeqn

which is not quite obvious from the construction of $O$. And so we get around
this problem by defining the algebra of {\it observables}, which morally
should be the smallest $k$-algebra extending the action of $A$ on
$\mathcal{C}$ and having all the nice functorial properties. We just claim
that the following definition works, which bluntly tries to fix the problem.

\begin{defn}  (algebra of observables)

The $k$-algebra of observables of the finite swarm $\mathcal{C}$ is the
sub-algebra 

\beqn
\sheaf (\mathcal{C},\pi ) := \underset{\mathcal{C}_0\subset\mathcal{C}}{\bigcap}
o(\mathcal{C}_0\subset\mathcal{C})^{-1} (O(\mathcal{C}_0,\pi )) \subset O
(\mathcal{C},\pi )
\eeqn

where $\mathcal{C}_0\subset\mathcal{C}$ runs through all sub-diagrams of
$\mathcal{C}$. Here $o(\mathcal{C}_0\subset\mathcal{C})^{-1}
(O(\mathcal{C}_0,\pi ))$ denotes the set-theoretic inverse image.
\end{defn}

Now following \cite{Lau3} we claim that this new $\sheaf$-construction of the
algebra of observables is a contravariant functor on the partially ordered
category of sub-diagrams of a given diagram $\mathcal{C}$. 

We need to extend the construnction to infinte diagrams. We call a swarm $\mathcal{C}$ {\it
  permissible} if there exists a $k$-algebra homomorphism,

\beqn
\eta(|\mathcal{C}|,\pi ): A\map O (|\mathcal{C}|,\pi ),
\eeqn
 
compatible with the morphisms $\eta (|\mathcal{C}_0|,\pi)$ and
$o(\mathcal{C}_0\subset \mathcal{C})$, where $\mathcal{C}_0$ runs through all
sub-diagrams of $\mathcal{C}$.

Now for a permissible (possibly infinite) swarm $\mathcal{C}$ we define

\beqn
\sheaf(\mathcal{C},\pi) :=
\underset{\mathcal{C}_0\subset\mathcal{C}}{\underset{\leftarrow}{Lim}} \,
\sheaf(\mathcal{C}_0,\pi ) 
\eeqn

where $\mathcal{C}_0$ runs through all finite sub-diagrams of $\mathcal{C}$.

Finally we arrive at the notion of the {\it structure sheaf} $\sheaf_\pi$. For
every finite sub-diagram $\mathcal{C}_0$ of a swarm $\mathcal{C}$, consider
the natural morphism 

\beqn 
\kappa (\mathcal{C}_0): \sheaf(\mathcal{C},\pi)\map (H(\mathcal{C}_0)\otimes
End(\mathcal{C})) \hspace{2mm} \text{refer to [\ref{Burnside}] and [\ref{pre-obs}]}
\eeqn

and consider the two-sided ideal $\mathtt{n}\subset \sheaf(\mathcal{C},\pi)$,
defined by 

\beqn
\mathtt{n} = \underset{\mathcal{C}_0\subset\mathcal{C}}{\bigcap} ker
\kappa(\mathcal{C}_0)
\eeqn

where $\mathcal{C}_0$ once again runs through all finite sub-diagrams of
$\mathcal{C}$. We now set

\beqn
\sheaf_\pi (\mathcal{C}) = \sheaf(\mathcal{C},\pi)/\mathtt{n}
\eeqn

Going through all these rigmarole should have some tangible benefit and that is given by
the following definition of an {\it affine \nc scheme}.

\begin{defn} (Affine non-commutative scheme)

A swarm of $A$-modules, $\mathcal{C}$, is called an affine scheme for $A$, if 

\beqn
\eta(\mathcal{C}):A\map \sheaf_\pi (\mathcal{C})
\eeqn

is an isomorphism and we shall refer to $A$ as the affine ring of this scheme.

\noindent
Note that we are talking about ``an'' affine scheme for $A$ and not ``the''
affine scheme as there could be several different scheme structures for $A$.
\end{defn}
 
\begin{rem}
For a commutative $k$-algebra $A$ (not necessarily finite type) let $Simp(A)$ denote the set of simple $A$-modules. It is seen that $\sheaf (Simp(A))\simeq \underset{\mathtt{m}}{\prod} \hat{A}_{\mathtt{m}}$, where $\mathtt{m}$ runs through all maximal ideals of $A$. However, $\eta (Simp(A))$ fails to be an isomorphism between $A\map \sheaf (Simp(A),\pi )$ in general. This can be rectified by just including $A$ in $Simp(A)$, which we denote by $Simp^\ast (A)$. 
\end{rem}

\spc
\begin{center}
\underline{\bf Infinitesimal Structures}
\end{center}

\vspace{2mm}

Let $\mathcal{C}$ be a swarm of $A$-modules and consider the point $x = V_i \in \mathcal{C}$. Set $\mathtt{V} = \underset{V\in \mathcal{C}}{\oplus} V$. 

\begin{defn} (Tangent space)

The tangent space at $x$, denoted by $T_x$ is defined to be

\beqn
\begin{split}
T_x \, :=\, & \{ \xi \in Ext^1_A (V_i,V_i)\; | \;\text{$\forall p$ $\exists \xi_p$ such that} \\
       & \text{$\forall \phi :=\phi_{i,p}:V_i\map V_p$ $\phi_\ast (\xi) = \phi^\ast (\xi_p)$ and} \\
       & \text{$\forall \phi :=\phi_{p,i}:V_p\map V_i$ $\phi_\ast (\xi_p) = \phi^\ast (\xi)$}\}
\end{split}
\eeqn

\end{defn}   

 It can be verified that the composition of the natural maps $Der_k (A,A)\map Der_k (A,End_k (V_i))$ and the surjection $Der_k (A,End_k (V_i))\map Ext^1_A (V_i,V_i)$ actually gives a map 

\beqn 
\theta_x : Der_k (A,A) \map T_x \subset Ext^1_A (V_i,V_i)
\eeqn

\begin{defn}  (Smoothness)

We call a point $x = V_i \in \mathcal{C}$ smooth if the map $\theta_x$ is surjective. If this is true for all points of $\mathcal{C}$, we say that $\mathcal{C}$ is smooth. 
\end{defn}

The author also has a short treatment of invariant theory and moduli at the end of \cite{Lau3}. All this is not futile but, due to lack of time, we need to end abruptly. The most interesting scheme structure for an algebra $A$ is probably given by $Simp^\ast(A)$, which has been discussed at length in \cite{Lau3}. These seemingly preternatural definitions have some interesting and revealing applications which can be found in \cite{JLS}. 


\newpage

\section{Synopsis of some other points of view}

\vspace{5mm}

This section is really meant to be a glossary rather than an in-depth
 analysis. In spite of that, given the staggering number of articles in this
 area, we shall only be able to see the tip of the iceberg. It is worth
 mentioning that curious readers may take a look at the web-page of {\it Paul
 J. Smith}

\begin{center} 
$http://www.math.washington.edu/~smith/Research/research.html$
\end{center}

\noindent
where one can find links to a host of other interesting web-sites catering to research in \nc algebra and algebraic geometry. 

\vspace{2mm}
\noindent
{\bf 1.} The first attempt that I would like to talk about is {\it Noncommutative
  spaces and flat descent} by {\it Kontsevich}
and {\it Rosenberg} and it uses absolutely ``state-of-the-art machinery''
(refer to \cite{KR2}). In a previous article entitled ``Noncommutative smooth spaces'' (see \cite{KR1}) the same authors had advocated two principles in \nc geometry. One principle asserts that \nc notions should come from the representations of their commutative counterparts and the other says that it should be possible to read off a \nc object from its ``covering datum''. The article under consideration improves upon the second line of thought, but it is much more abstract in nature. This approach adheres to sheaves of sets on the category of affine schemes endowed with flat topology and having descent property. The minimum desideratum for \nc geometry, or for
that matter any geometric theory, is a category of local or ``affine''
objects, and a functor which to every such object assigns a ``space'', like
$Spec$. This article relies on the philosophy that a ``space'' can
be identified with categories, which are thought of as categories of coherent
or quasi-coherent sheaves on it. So the authors carry along with each ``space'' $X$
a category $\mathcal{C}_X$, which is regarded as the category of quasi-coherent sheaves
on $X$. So it is a mathematical object like a fibred category

\[
\begin{CD}
\mathcal{C} \\
@V\pi VV \\
\mathcal{B}
\end{CD}
\]

where the base category $\mathcal{B}$ serves as the category of  
``spaces''. For any $X\in\mathcal{B}$  we denote by $\mathcal{C}_X$ the fibre
category over $X$, which is thought of as the category of quasi-coherent
sheaves on $X$. For any morphism $X\overset{f}{\map} Y$ we denote by $f^\ast$
the corresponding $2$-morphism $\mathcal{C}_Y\map\mathcal{C}_X$ and call it
the {\it inverse image functor of $f$}. 

Now we set $\mathcal{C} = Cat^{op}$ and $\mathcal{B} = |Cat|^0$ where $Cat$
stands for the $2$-category of categories and we set about defining
$|Cat|^0$. 

\vspace{1mm}
\noindent
$Ob(|Cat|^0) = $ objects of $Cat$ {\it i.e.,} categories. 

\noindent
$Hom_{|Cat|^0}(X,Y) = $ isomorphism classes of functors from $Cat^{op}_Y\map
Cat^{op}_X$. 

\vspace{1mm}
The authors write $f=[F]$ to indicate that $f$ is a morphism such
that the functor $F$ belongs to the isomorphism class of its inverse image
functors. Define a functor from $Cat^{op}\map |Cat|^0$, which is identity on
the objects and which sends each morphism $Cat^{op}_Y\overset{F}{\map}
Cat^{op}_X$ to the morphism 
$X\overset{[F]}{\map} Y$, to make $Cat^{op}$ into a fibred category over
$|Cat|^0$. The authors then orchestrate a definition of a ``cover''. Using
these covers they are able to define affine objects and gluing of such objects to form locally affine objects, which can be
recovered uniquely from the covering data via ``flat descent''. To satisfy one's curiosities, interested readers are encouraged to go through this article.

\vspace{2mm}
\noindent
{\bf 2.} Now we are going to talk about {\it Noncommutative Geometry Based on
  Commutator Expansions} by {\it M. Kapranov} (refer to \cite{Kap}). The
  author himself describes it as an attempt to develop \nc algebraic geometry
  ``in the perturbative regime'' around ordinary commutative geometry. Only algebras over the complex numbers are considered. So given any associative $\CC$-algebra $R$, a decreasing filtration of $R$, denoted by $\{F^dR\}_{d\geqslant
  0}$ is defined as follows:

\beqn
F^dR=\underset{m}{\sum}\underset{(i_1 + \dots + i_m - m = d)}{\sum}
R.R^{Lie}_{i_1}.R\dots R.R^{Lie}_{i_m}.R
\eeqn

where $R^{Lie}$ is the Lie algebra concocted from the associative algebra
$R$. Then a suitable notion of $NC$-completion ($R\map R_{NC-complete} = \underset{\leftarrow}{Lim}\, R/F^d R$) of such algebras is defined
and for an $NC$-complete algebra $R$ he is able to define the localization for
any multiplicatively closed $T\subset R_{ab}$ ($R_{ab}:=R/[R,R]$) as:

\beqn
R\| T^{-1}\| = \underset{\leftarrow}{Lim}\, (R/F^{d+1}R)[T^{-1}]
\eeqn

For every $NC$-nilpotent algebra a space with a sheaf of functions
$X=\mathtt{Spec}\, (R) = (Spec\, R_{ab},\sheaf_X)$ can be defined as follows:

\beqn
\sheaf_X (D_g)= R[g^{-1}]
\eeqn 

where $D_g = \{ \wp\in X | g \notin \wp \}$. And the stalks are defined as 

\beqn
\sheaf_{X,\wp} := R_{\wp} = \underset{\underset{\wp\in D_g}{\rightarrow}}{Lim}
R[g^{-1}]
\eeqn

Set $X_{ab} = Spec\, R_{ab}$. Then, for an $NC$-complete algebra $R$, the formal spectrum $Spf(R)$ is given by the ringed space $(X_{ab},\sheaf_X)$, where $\sheaf_X$ is the sheaf of
topological rings obtained as the inverse limit of the structure sheaves of
$\mathtt{Spec}\, (R/F^{d+1}R)$.

The ringed space $Spf(R)$ is called an affine $NC$-scheme and $\mathcal{NC}$
is the category of $NC$-complete algebras and, as expected, there is indeed an
equivalence of categories between the category of $NC$-schemes and
$\mathcal{NC}$ given by $X\map\Gamma (X_{ab},\sheaf_X)$. The story does not
end here. Given an ordinary manifold $M$, he considers objects (thickenings)
$X$ with $X_{ab} = M$ and towards the end illustrates that several familiar
algebraic varieties, including the classical flag varieties and all the smooth
moduli spaces of vector bundles, possess natural $NC$-thickenings. In addition, notions like $NC$-smoothness are discussed and, conforming to our expectations, for an $NC$-smooth algebra $R$, $dim\, R=n$, if $x$ is a $\CC$-valued point of $Spec\, R_{ab}$, it is shown that $\hat{R}_x \simeq \CC\langle\langle x_1,\dots ,x_n\rangle\rangle$.

\vspace{2mm}
\noindent
{\bf 3.} An book published in 1981 entitled ``Noncommutative Algebraic Geometry'' by {\it Oystaeyen} and {\it Verschoren} (refer to \cite{OV}) promulgates a theory which works well for affine $PI$-algebras over an algebraically closed field $k$. It heavily relies on the concomitant developments in \nc ring theory, carried out by {\it Artin}, {\it Schelter}, {\it Cohn}, {\it Procesi} and the authors themselves, to name only a few. One of the watermarks of the book is the \nc version of {\it Riemann-Roch} for curves.  

\vspace{2mm}
\noindent
{\bf 4.} We may ``derive'' the philosophy of identifying a space with its category of quasi-coherent sheaves. Ever since the article \cite{Bei} by {\it Be\u\i linson} describing the derived category of coherent sheaves on $\proj^n$ there has been a spate of activities on such things. We would like to point out particularly the article entitled ``Derived categories of coherent sheaves'' by {\it Bondal} and {\it Orlov} which ends with some intriguing connections with \nc geometry (refer to \cite{BO3}). They try to do some kind of birational geometry in this derived set-up and {\it Orlov} has taken this notion one step further with his article \cite{Orl1}. Grothendieck's notions of a ``site'' and a ``topos'' hold the promise of being generalized to the \nc framework and it seems that something similar has been achieved by {\it Orlov} in \cite{Orl2}. 

\vspace{2mm}
\noindent
{\bf 5.} A highly readable article by {\it Pierre Cartier} \cite{Car} entitled ``A MAD DAY'S WORK: FROM GROTHENDIECK TO CONNES AND KONTSEVICH. THE EVOLUTION OF CONCEPTS OF SPACE AND SYMMETRY'' provides a sublime description of the evolution of our notions of space and symmetry. It also discusses some connections between \nc geometry via topos theory after {\it Grothendieck} and that of {\it Connes} using operator algebras (He refers to \cite{Tap} in the process). Interestingly, it begins with a short biography of {\it Alexander Grothendieck} and end with ``a dream'' of {\it Cartier}. It cannot be totally classified as an article in \nc algebraic geometry but the whole treatise is not only instructive but edifying. 

\vspace{2mm}
Finally we just mention some other sources of knowledge in this area. There is a set of notes by {\it Le Bruyn} \cite{Bru} which can be downloaded freely from his web-site. {\it Smith} and {\it Zhang} had disseminated a nice approach in \cite{SZ} of studying curves in \nc schemes. Lately {\it Polishchuk} has come up with some interesting ideas for which two good sources could be \cite{Pol1} and \cite{Pol2}. {\it M. Van den Bergh} has an article entitled ``Blowing up of non-commutative smooth surfaces'' \cite{VDB1} and he has followed it up with a recent article in the arXiv server, dealing with blowing down \cite{VDB2}. Incidentally {\it Keeler}, {\it Rogalski} and {\it Stafford} also have their own version of blowing up in \cite{KRS}. Last and perhaps the closest to being definitive article providing an overview of the goings-on in this area is \cite{StV} by {\it Stafford} and {\it Van den Bergh}, which was also the backbone of this write-up. Recently {\it Matilde Marcolli} has written a book entitled ``Arithmetic Noncommutative Geometry'' \cite{Mar}, which explores some interesting and intriguing connections between \nc geometry and arithmetic. 

There are also some examples of \nc spaces like $D$-schemes of {\it Be\u\i linson-Bernstein} and \nc schemes of {\it P. Cohn} which readers are invited to unravel for themselves.

\vspace{2mm}
Noncommutative geometry is now a rich and widely pursued subject. If nature desires to divest herself of commutativity, as examples coming from physics manifest, then it may not be as ``natural'' as we think. One question that was mulled by several people is whether noncommutative geometry is just what we obtain once we purge the commutativity hypothesis. Also, is \nc geometry just a shadow of its more resplendent predecessor? There is room for arbitration but clearly, it has had some impact on our perception of a ``space'' and it will continue to engage us in its myriad problems.


\newpage

\section{ACKNOWLEDGEMENTS}

\vspace{5mm}

So far a meticulous reader must have noticed that the first person singular number was very carefully avoided. However, here I would just be myself and do the honours. 

First of all I would like to thank my advisor, {\it Prof. Matilde Marcolli},
for introducing me to this subject and lending me a patient ear during the discussion sessions that we had. She was also instrumental in
helping me attend a few schools and conferences, which have substantially
augmented my knowledge, and edited my write-up on quite a few occasions. 

I would also like to thank some of my friends, particularly {\it Jorge Andres
  Plazas Vargas} and {\it Eugene Ha} in Max-Planck-Institut f\" ur Mathematik,
Bonn who took their time out of their busy schedules to take part in our
regular discussion sessions and often providing insightful remarks.

A would specially like to thank {\it Prof. Michel Van den Bergh} for going
through the first two chapters and providing valuable comments, suggestions
and sometimes even graphic explanations to my queries. Even while writing the later chapters we communicated via e-mail and he was always forthcoming with his answers, however silly my questions might have been. 

It would be thoroughly unfair if I fail to acknowledge the people who have
proof read this article at various stages of its development or helped me in LaTeX-ing, namely {\it Gunther Vogel} and {\it \" Ozg\" ur Ceyhan}. Gunther had also provided me with some useful references.  

Finally, I take this opportunity to extend my earnest gratitude to everyone at
Max-Planck-Institut f\" ur Mathematik, Bonn, for providing a very conducive
environment for studies.


\bibliographystyle{alpha}
\bibliography{ncgSurvey}

\begin{thebibliography}{AVdB90}

\bibitem[AVdB90]{AV}
M.~Artin and M.~Van~den Bergh.
\newblock Twisted homogeneous coordinate rings.
\newblock {\em J. Algebra 133, 249-271}, MR 91k:14003, 1990.

\bibitem[AZ94]{AZ}
M.~Artin and J.~J. Zhang.
\newblock Noncommutative projective schemes.
\newblock {\em Adv. in Math. 109, no. 2, 228-287}, MR 96a:14004, 1994.

\bibitem[Be{\u\i}84]{Bei}
A.~A. Be{\u\i}linson.
\newblock The derived category of coherent sheaves on \text{$\proj^n$}.
\newblock {\em Selected translations. Selecta Math. Soviet, 3, no. 3, 233-237},
  MR 88h:14021, 1983/84.

\bibitem[BK89]{BK}
A.~I. Bondal and M.~M. Kapranov.
\newblock Representable functors, serre functors and reconstructions.
\newblock {\em Izv. Akad. Nauk SSSR Ser. Mat. 53, no. 6, 1183-1205, 1337}, MR
  91b:14013, 1989.

\bibitem[BO01]{BO1}
A.~I. Bondal and D.~Orlov.
\newblock Reconstruction of a variety from the derived category and groups of
  autoequivalences.
\newblock {\em Compositio Math. 125, no. 3, 327-344}, MR 2001m:18014, 2001.

\bibitem[BO02]{BO3}
A.~Bondal and D.~Orlov.
\newblock Derived categories of coherent sheaves.
\newblock {\em Proceedings of the International Congress of Mathematicians,
  Vol. 2, 47-56, Higher Ed. Press, Beijing}, MR 2003m:18015, 2002.

\bibitem[Bruyn]{Bru}
L.~Le Bruyn.
\newblock Noncommutative geometry@n.
\newblock {\em "Forgotten" book, available from the author's homepage},
  \text{{\it http://www.math.ua.ac.be/$\sim$lebruyn/}}.

\bibitem[BVdB03]{BV}
A.~I. Bondal and M.~Van~den Bergh.
\newblock Generators and representability of functors in commutative and
  noncommutative geometry.
\newblock {\em Moscow Math. J. vol. 3, no. 1, 1-36, 258}, 2003.

\bibitem[Car01]{Car}
Pierre Cartier.
\newblock A mad day's work: from \text{Grothendieck} to \text{Connes} and
  \text{Kontsevich}.
\newblock {\em Bull. Amer. Math. Soc.(N.S.) 38, no. 4, 389-408}, MR
  2000c:01028, 2001.

\bibitem[Con94]{Con}
A.~Connes.
\newblock {\em Noncommutative Geometry}.
\newblock Academic Press, Inc., San Diego, CA., MR 95j:46063, 1994.

\bibitem[Gab62]{Gab}
P.~Gabriel.
\newblock Des cat\' egories ab\' eliennes.
\newblock {\em Bull. Soc. Math. France 90, 323-448}, MR 38:1411, 1962.

\bibitem[Gro57]{Gro}
A.~Grothendieck.
\newblock Sur quelques points d'alg\` ebre homologique.
\newblock {\em T{\^{o}}huku Math. J. 9, 119-221}, MR 21:1328, 1957.

\bibitem[Har77]{Har}
R.~Hartshorne.
\newblock {\em Algebraic Geometry}.
\newblock Springer-Verlag, 1977.

\bibitem[Ile90]{Ile}
R.~Ile.
\newblock Noncommutative deformations. \text{$ab\neq ba$}.
\newblock {\em Master thesis, \text{University} of \text{Oslo}}, 1990.

\bibitem[J{\o}r97]{Jor1}
P.~J{\o}rgensen.
\newblock Serre duality for \text{Tails(A)}.
\newblock {\em Proc. Amer. Math. Soc. 125, 709-716}, MR 97e:14002, 1997.

\bibitem[J{\o}r00]{Jor3}
P.~J{\o}rgensen.
\newblock Intersection theory on noncommutative surfaces.
\newblock {\em Trans. Amer. Math. Soc., 352, 5817-5854}, CMP 99:14, 2000.

\bibitem[J{\o}r10]{Jor2}
P.~J{\o}rgensen.
\newblock Noncommutative graded homological identities.
\newblock {\em J. London Math. Soc., 57, 336-350}, MR 99h:16010.

\bibitem[JZ00]{JZ}
P.~J{\o}rgensen and J.~J. Zhang.
\newblock Gourmet guide to \text{Gorensteinness}.
\newblock {\em Adv. Math., 151, 313-345}, CMP 2000:12, 2000.

\bibitem[Kap98]{Kap}
M.~Kapranov.
\newblock Non-commutative geometry based on commutator expansions.
\newblock {\em J. Reine Angew. Math. 505, 73-118}, MR 2000b:14003, 1998.

\bibitem[KL85]{KL}
G.~R. Krause and T.~H. Lenagan.
\newblock Growth of algebras and \text{Gelfand-Kirillov} dimension.
\newblock {\em Research Notes in Math., vol 116, Pitman, Boston}, MR 86g:16001,
  1985.

\bibitem[KR00]{KR1}
M.~Kontsevich and A.~Rosenberg.
\newblock Noncommutative smooth spaces.
\newblock {\em The Gelfand Mathematical Seminars, 1996-1999, 85-108, Gelfand
  math. Sem., Birkh\" auser, Boston, MA}, CMP 2000:17, 2000.

\bibitem[KR04]{KR2}
M.~Kontsevich and A.~Rosenberg.
\newblock Noncommutative spaces and flat descent.
\newblock {\em MPIM Preprint}, 2004.

\bibitem[Lau79]{Lau1}
O.~A. Laudal.
\newblock Formal moduli of algebraic structures.
\newblock {\em Lecture Notes in Math., 754, Springer, Berlin}, MR 82h:14009,
  1979.

\bibitem[Lau00]{Lau3}
O.~A. Laudal.
\newblock Noncommutative algebraic geometry.
\newblock {\em MPIM Preprint, 115}, 2000.

\bibitem[Lau02]{Lau2}
O.~A. Laudal.
\newblock Noncommutative deformations of modules.
\newblock {\em The Roos Festschrift, vol. 2, Homology Homotopy Appl. 4, no. 2,
  part 2, 357-396}, MR 2003e:16005, 2002.

\bibitem[LS04]{JLS}
S.~J{\o}ndrup O.~A. Laudal and A.~B. Sletsj{\o}e.
\newblock Noncommutative plane curves.
\newblock {\em Institute Mittag-Leffler, Report No. 20}, 2003/2004.

\bibitem[Man88]{Man2}
Yu.~I. Manin.
\newblock Quantum groups and noncommutative geometry.
\newblock {\em Univ. de Montr\' eal, Centre de Recherches Math.,Montreal, QC},
  MR 91e:17001, 1988.

\bibitem[Man91]{Man1}
Y.~I. Manin.
\newblock {\em Topics in Noncommutative Geometry}.
\newblock Princeton University Press, New Jersey, MR 92k:58024, 1991.

\bibitem[Mar05]{Mar}
Matilde Marcolli.
\newblock {\em Arithmetic noncommutative geometry}, volume~36 of {\em
  University Lecture Series}.
\newblock American Mathematical Society, Providence, RI, 2005.
\newblock With a foreword by Yuri Manin.

\bibitem[Orl03a]{Orl1}
D.~Orlov.
\newblock Derived categories of coherent sheaves and equivalences between them
  (russian).
\newblock {\em Uspekhi Mat. Nauk, 58, no. 3, 89-172; translation in Russian
  Math. Surveys, 58, no. 3, 511-591}, MR 2004g:14021, 2003.

\bibitem[Orl03b]{Orl2}
D.~Orlov.
\newblock Quasicoherent sheaves in commutative and noncommutative geometry
  (russian).
\newblock {\em Izv. Ross. Akad. Nauk Ser. Mat., 67, no. 3, 119-138; translation
  in Izv. Math., 67, no. 3, 535-554}, MR 2004m:14026, 2003.

\bibitem[OV81]{OV}
F.~Van Oystaeyen and A.~Verschoren.
\newblock {\em Noncommutative Algebraic Geometry}.
\newblock Lecture Notes in Math., vol 887, Springer Verlag, Berlin, MR
  85i:16006, 1981.

\bibitem[Pol04a]{Pol2}
A.~Polishchuk.
\newblock Noncommutative two-tori with real multiplication as noncommutative
  projective varieties.
\newblock {\em J. Geom. Phys., 50, no. 1-4, 162-187}, MR 2078224, 2004.

\bibitem[Pol04b]{Pol1}
A.~Polishchuk.
\newblock Classification of holomorphic vector bundles on noncommutative
  two-tori.
\newblock {\em Doc. Math., 9, 163-181}, MR2054986, 2004.

\bibitem[Rei75]{Rei}
I.~Reiner.
\newblock {\em Maximal Orders}.
\newblock Academic Press, New York, MR 52:13910, 1975.

\bibitem[Rin84]{Rin}
C.~M. Ringel.
\newblock {\em Tame algebras and integral quadratic forms}.
\newblock Lecture Notes in Math., vol. 1099, Springer Verlag, Berlin, MR
  87f:16027, 1984.

\bibitem[Rog02]{Rog}
D.~Rogalski.
\newblock {\em Examples of generic noncommutative surfaces}.
\newblock PhD thesis, \text{University of Michigan}, MR 2004b:58007, May 2002.

\bibitem[Ros95]{Ros1}
A.~Rosenberg.
\newblock {\em Noncommutative Algebraic Geometry and Representations of
  Quantized Algebras}.
\newblock Mathematics and its Applications, vol. 330, Kluwer Academic
  Publishers, Dordrecht, MR 97b:14004, 1995.

\bibitem[Ros98]{Ros2}
A.~Rosenberg.
\newblock The spectrum of abelian categories and reconstructions of schemes.
\newblock {\em Rings, Hopf Algebras, and Brauer groups, Lectures Notes in Pure
  and Appl. Math., vol. 197, Marcel Dekker, New York, 257-274}, MR 99d:18011,
  1998.

\bibitem[RS03]{KRS}
D.~S. Keeler~D. Rogalski and J.~T. Stafford.
\newblock Na\" ive noncommutative blowing up.
\newblock {\em arXiv:math.RA/0306244 v1}, 2003.

\bibitem[RVdB02]{RV}
I.~Reiten and M.~Van~den Bergh.
\newblock Noetherian hereditary categories satisfying serre duality.
\newblock {\em J. Amer. Math. Soc. 15, no. 2, 295-366}, MR 2003a:18011, 2002.

\bibitem[Sch68]{Sch}
M.~Schlessinger.
\newblock Functors of artin rings.
\newblock {\em Trans. Amer. Math. Soc. 130, 208-222}, MR0217093, 1968.

\bibitem[Ser55]{Ser}
J.~P. Serre.
\newblock Faisceaux alg\' ebriques coh\' erents.
\newblock {\em Ann. of Math. (2) 61, 197-278}, MR 16:953c, 1955.

\bibitem[Ste75]{Ste}
B.~Stenstr{\"o}m.
\newblock Rings of quotients.
\newblock {\em Die Grundlehren der mathematischen Wissenschaften in
  Einzeldarstellungen, vol. 217, Springer Verlag, Berlin}, MR 52:10782, 1975.

\bibitem[SV00]{FSV}
E.~M. Friedlander~A. Suslin and V.~Voevodsky.
\newblock Cycles, transfers, and motivic homology theories.
\newblock {\em Ann. of Math. Stud., 143, Princeton Univ. Press, Princeton, NJ},
  MR 2001d:14026, 2000.

\bibitem[SVdB01]{StV}
J.~T. Stafford and M.~Van~den Bergh.
\newblock Noncommutative curves and noncommutative surfaces.
\newblock {\em Bull. Amer. Math. Soc. (N.S.) 38, no. 2, 171-216}, MR
  2002d:16036, 2001.

\bibitem[SZ94a]{StZ1}
J.~T. Stafford and J.~J. Zhang.
\newblock Examples in non-commutative projective geometry.
\newblock {\em Math. Proc. Cambridge Phil. Soc. 116, 415-433}, MR 95h:14001,
  1994.

\bibitem[SZ94b]{StZ2}
J.~T. Stafford and J.~J. Zhang.
\newblock Homological properties of (graded) \text{Noetherian} \text{PI} rings.
\newblock {\em J. Algebra, 168, 988-1026}, MR 95h:16030, 1994.

\bibitem[SZ98]{SZ}
S.~P. Smith and J.~J. Zhang.
\newblock Curves on quasi-schemes.
\newblock {\em Alg. Rep. Theory 1, 311-351}, CMP 99:11, 1998.

\bibitem[Tap91]{Tap}
J.~Tapia.
\newblock Quelques spectres en \text{K}-th\' eorie topologique des alg\` ebres
  de fr\' echet et applications \` a l'alg\` ebre des fonctions de classe
  \text{$C^\infty$} sur une vari\' et\' e.
\newblock {\em Pr\' eprint IHES}, IHES/M/91/37, 1991.

\bibitem[TVdB90]{ATV1}
M.~Artin~J. Tate and M.~Van~den Bergh.
\newblock Some algebras associated to automorphisms of elliptic curves.
\newblock {\em The Grothendieck Festschrift, vol. 1, Birkh\" auser, Boston,
  33-85}, MR 92e:14002, 1990.

\bibitem[TVdB91]{ATV2}
M.~Artin~J. Tate and M.~Van~den Bergh.
\newblock Modules over regular algebras of dimension $3$.
\newblock {\em Invent. Math. 106, 335-388}, MR 93e:16055, 1991.

\bibitem[VdB98]{VDB2}
M.~Van~den Bergh.
\newblock Abstract blowing down.
\newblock {\em arXiv:math.AG/9806046 v1}, 1998.

\bibitem[VdB01]{VDB1}
M.~Van~den Bergh.
\newblock Blowing up of non-commutative smooth surfaces.
\newblock {\em Mem. Amer. Math. Soc., 154, no. 734}, MR 2002k:16057, 2001.

\bibitem[Ver92]{Ver}
A.~B. Verevkin.
\newblock On a non-commutative analogue of the category of coherent sheaves on
  a projective scheme.
\newblock {\em Amer. Math. Soc. Transl. 151, 41-53}, MR 93j:14002, 1992.

\bibitem[WZ00]{WZ}
Q.~S. Wu and J.~J. Zhang.
\newblock Some homological invariants of local \text{PI} algebras.
\newblock {\em J. Algebra 225, 904-935}, CMP 2000:09, 2000.

\bibitem[YZ97a]{YZ2}
A.~Yekutieli and J.~J. Zhang.
\newblock Rings with auslander dualizing complexes.
\newblock {\em J. Algebra 213, 1-51}, MR 2000f:16012, 1997.

\bibitem[YZ97b]{YZ1}
A.~Yekutieli and J.~J. Zhang.
\newblock Serre duality for noncommutative projective schemes.
\newblock {\em Proc. Amer. Math. Soc. 125, 697-707}, MR 97e:14003, 1997.

\bibitem[Zha96]{Zha}
J.~J. Zhang.
\newblock Twisted graded algebras and equivalences of graded categories.
\newblock {\em Proc. London Math. Soc. (3), 72, no. 2, 281-311}, MR 96k:16078,
  1996.

\end{thebibliography}

\end{document}